%Paper: q-alg/9512027
%From: Atsushi Nakayashiki <atsushi@rcmath.rc.kyushu-u.ac.jp>
%Date: Sun, 24 Dec 1995 21:27:29 +0900

%
%--------------------------------------------------------macros
%
%
%   version 28 June 1994
%

\def\mmmth{\mathsurround=0pt}

\def\fsquare(#1,#2){
\hbox{\vrule$\hskip-0.4pt\vcenter to #1{\normalbaselines\mmmth
\hrule\vfil\hbox to #1{\hfill$\scriptstyle #2$\hfill}\vfil\hrule}$\hskip-0.4pt
\vrule}}

\def\addsquare(#1,#2){\hbox{$
	\dimen1=#1 \advance\dimen1 by -0.8pt
	\vcenter to #1{\hrule height0.4pt depth0.0pt%
	\hbox to #1{%
	\vbox to \dimen1{\vss%
	\hbox to \dimen1{\hss$\scriptstyle~#2~$\hss}%
	\vss}%
	\vrule width0.4pt}%
	\hrule height0.4pt depth0.0pt}$}}

\def\Addsquare(#1,#2){\hbox{$
	\dimen1=#1 \advance\dimen1 by -0.8pt
	\vcenter to #1{\hrule height0.4pt depth0.0pt%
	\hbox to #1{%
	\vbox to \dimen1{\vss%
	\hbox to \dimen1{\hss$~#2~$\hss}%
	\vss}%
	\vrule width0.4pt}%
	\hrule height0.4pt depth0.0pt}$}}

\def\Fsquare(#1,#2){
\hbox{\vrule$\hskip-0.4pt\vcenter to #1{\normalbaselines\mmmth
\hrule\vfil\hbox to #1{\hfill$#2$\hfill}\vfil\hrule}$\hskip-0.4pt
\vrule}}

\def\naga{%
	\hbox{$\vcenter to 0.4cm{\normalbaselines\mmmth
	\hrule\vfil\hbox to 1.2cm{\hfill$\cdots$\hfill}\vfil\hrule}$}}

\def\Flect(#1,#2,#3){
\hbox{\vrule$\hskip-0.4pt\vcenter to #1{\normalbaselines\mmmth
\hrule\vfil\hbox to #2{\hfill$#3$\hfill}\vfil\hrule}$\hskip-0.4pt
\vrule}}

\def\PFlect(#1,#2,#3){
\hbox{$\hskip-0.4pt\vcenter to #1{\normalbaselines\mmmth
\vfil\hbox to #2{\hfill$#3$\hfill}\vfil}$\hskip-0.4pt}}

\dimen1=0.5cm\advance\dimen1 by -0.8pt

\def\vnaka{\normalbaselines\mmmth\baselineskip0pt\offinterlineskip%
	\vrule\vbox to 0.6cm{\vskip0.5pt\hbox to
\dimen1{$\hfil\vdots\hfil$}\vfil}\vrule}

\dimen2=1.5cm\advance\dimen2 by -0.8pt

\def\vnakal{\normalbaselines\mmmth\baselineskip0pt\offinterlineskip%
	\vrule\vbox to 1.2cm{\vskip7pt\hbox to
\dimen2{$\hfil\vdots\hfil$}\vfil}\vrule}

\hsize=13cm
\magnification=\magstep1
%
%%%%%%%%%%%%%%%%%%%%%%%%%%%%%%%%%%%%%%%%%%%%%%%%%%%%%%%%%%%%%%%%%

\message{FrameTeX Ver 1.0 by M. Hashimoto}
%
%  Usage.... \Young<size of a box>(partition)    (see an example below)
%            \young<size of a box>(partition)    (for subscript)
%
%\aho       ... temporary macro
%\b@xsize   ... absolute size of one box in the Young diagram, default=6pt.
%               default=2.5pt for \young
%\p@rtition ... a sequence of nonnegative integers, separated by commas
%\getnum    ... get a number of boxes in the current row and put it into
%               \r@. it removes this number and the following comma (if any).
%               it must be \global.
%\if@cont   ... If \p@rtition is exhausted by \getnum it gets to be
%               false it must be \global.
%\r@        ... The real value of \r@ is 2*(the number of boxes).
%               see \getnum. it must be \global.
%\Young     ... The macro that creates a Young diagram
%\@@hline   ... Creates a horizontal line.
%\r@wofbox  ... Creates a row of boxes.
%\@nd       ... Equivalent to &. But this token list is not forbidden.
%\@nds      ... Token list consisting of \@nd's.
\catcode`\@11 \newdimen\b@xsize
\newtoks\p@rtition
\newcount\r@
\newtoks\@nds
\newif\if@cont
\def\Young#1(#2){\def\aho{#1}\ifx\aho\empty\b@xsize6pt\else\b@xsize\aho\fi
\global\p@rtition{#2,\aho}\global\@conttrue\expandafter
\getnum\the\p@rtition
\ifnum\r@>0\global\setbox1\hbox to\b@xsize{\hfil
\vrule width0ptheight\b@xsize depth0pt}
\vcenter{\hbox{\vbox{\offinterlineskip\tabskip0pt\def\@nd{&}
\halign{\vrule height\b@xsize##&&\copy1##&\vrule##\cr\@@hline\thisline@\cr}}}}
\else\emptyset\fi}%
\def\young#1(#2){\def\aho{#1}\ifx\aho\empty\b@xsize2.5pt\else\b@xsize\aho\fi
\global\p@rtition{#2,\aho}\global\@conttrue\expandafter
\getnum\the\p@rtition
\ifnum\r@>0\global\setbox1\hbox to\b@xsize{\hfil
\vrule width0ptheight\b@xsize depth0pt}\hbox{\vtop{\offinterlineskip\tabskip0pt
\def\@nd{&}\halign{\vrule height\b@xsize##&&\copy1##&\vrule##\cr
\@@hline\thisline@\cr}}}\else\emptyset\fi}%
\def\thisline@{\cr\r@wofbox\cr\@@hline\if@cont\let\noriP
\tr@ns\else\let\noriP\relax\fi\noriP}%
\def\tr@ns{\expandafter\getnum\the\p@rtition\thisline@}%
\def\getnum#1,#2\aho{\def\aho{#2}\ifx\aho\empty\def\aho{0}\global\@contfalse\fi
\ifnum#1=0\global\@contfalse\fi
\global\r@#1\global\multiply\r@\tw@
\global\p@rtition{#2\aho}}%
\def\@@hline{\omit\mscount\r@\loop\ifnum\mscount>\@ne\sp@n\repeat\hrulefill&}
\def\r@wofbox{\global\@nds{}\global\count\@ne\r@\gdef\foolish{\global\@nds
\expandafter{\the\@nds\@nd}\global\advance\count\@ne -\@ne
\ifnum\count\@ne>0}\boke\the\@nds}%
%% FOLLOWING LINE CANNOT BE BROKEN BEFORE 80 CHAR
\def\boke{\foolish\global\let\ahoaho\boke\else\global\let\ahoaho\relax\fi\ahoaho}%
\catcode`\@\active

%%%%%%%%%%%%%%%%%%% The end of frame.tex %%%%%%%%%%%%%%%

%\input young.tex
%\input ks.tex
%\magnification=\magstep1
\hsize=17truecm

\font\fontA=cmr12 scaled \magstep3
\font\fontB=cmr12 scaled \magstep2

\centerline{\fontA
Kostka Polynomials
}
\vskip1cm
\centerline{\fontA
and}
\vskip1cm
\centerline{\fontA
Energy Functions in Solvable Lattice Models
}
\vskip2cm
\centerline{
\rm
Atsushi Nakayashiki and Yasuhiko Yamada
}
\vskip1cm
\centerline{
\it
Graduate School of Mathematics, Kyushu University
}
\vskip2cm
\baselineskip 15pt
\noindent
{\bf Abstract}
\vskip2mm
\noindent
The relation between the charge of Lascoux-Schuzenberger
and the energy function in solvable lattice models is clarified.
As an application, A.N.Kirillov's conjecture on the expression of the
branching coefficient of ${\widehat {sl_n}}/{sl_n}$ as a
limit of Kostka polynomials is proved.
\vfill\eject

\pageno=1
%Nakayashiki personal macros
\def\ot{\otimes}
\def\ra{\longrightarrow}
\def\et{\tilde{e}}
\def\ft{\tilde{f}}
\def\La{\Lambda}
\def\la{\lambda}
\def\mapright#1{\smash{\mathop{\longrightarrow}\limits^{#1}}}
\def\ind{\hbox{ind}}

\noindent
{\fontB 1. Introduction}
\vskip2mm
\noindent
The study of two dimensional solvable
lattice models has provided many new links with other
areas of mathematics and physics.
In particular, one basic object of the models,
referred to as the one dimensional
configuration sum (1DCS) gives a generalization of
Rogers-Ramanujan type identities [ABF] and is related with
the characters of (quantum) affine Lie algebras [DJKMO].
The relation between 1DCS and affine Lie algebras was
systematized using the theory of crystal bases by Kashiwara [K1,KMN1-2].

The weight function describing the 1DCS is called the energy function
and it is determined from the $q=0$ behavior of the R matrix.
The theory of crystal bases gives a simple characterization
of energy functions without need for explicit forms
of R matrices [KMN1].
Obtaining this characterization played an essential role
in proving the relations between 1DCS and string functions [KMN1]
and the branching coefficients [DJO] of affine Lie algebras.

In this paper we add further application of the energy function
to combinatorics and representation theory.
In particular we obtain a new expression of the Kostka polynomials
in terms of energy functions. (for instance see (1.1) below).
This result was obtained while trying to generalize
our previous work [NY1-2] on spinon bases to the case of $\widehat{sl_n}$.

The Kostka polynomial $K_{\la\mu}(q)$ is the connection coefficient between
Schur function $s_\la(x)$ and the Hall-Littlewood symmetric function
$P_\mu(x:q)$:
$$
s_\la(x)=\sum_{\mu}K_{\la\mu}(q)P_\mu(x:q).
$$
The Kostka polynomial can be understood in many ways,
as the $q$ analogue of weight multiplicities [L], the $q$ analogue
of Clebsch-Gordan coefficients [Kir2], the Kazhdan-Lusztig polynomials,
[L] etc.
Among the ways of expressing the Kostka polynomial,
the combinatorial description of Lascoux-Schutzenberger [LS,B]
takes the form
$$
K_{\la\mu}(q)=\sum_{T\in{\cal T}(\la,\mu)}q^{c(T)},
$$
where ${\cal T}(\la,\mu)$ is the set of semi-standard tableaux
of shape $\la$ and weight $\mu$, the non-negative integer $c(T)$
is called the "charge", defined as a sum of the indices
(for a precise definition, see \S3 and [M] \S3).
Our main theorem (Theorem 4.1)
gives an expression of each of the indices in terms
of the energy function of $U_q(\widehat{sl_n})$.
Hence the charge of Lascoux-Schutzenberger appears to have
a representation theoretical meaning in quantum affine algebra.

To avoid notational complexity, in this introduction we
illustrate our results with one of their applications.
For the special case $\mu=(k^l)$, combining our results and
the theorem of Lascoux-Schutzenberger, we have
$$
K_{\la\mu}(q)=\sum_{(b_l,\cdots,b_1)\in{\cal T}(\la,\mu)}
q^{\sum_{j=1}^ljH_{k\La_1k\La_1}(b_{j+1},b_j)},
\eqno(1.1),
$$
where the sum is taken over $(b_l,\cdots,b_1)$, which is a path
(see \S2 for precise definition) parametrizing the set of
highest weight vectors with weight $\la$ in the tensor product
$V_{k\La_1}^{\ot l}$, $V_{k\La_1}$ being the irreducible representation
of $sl_n$ with highest weight $k\La_1$. The function
$H_{k\La_1k\La_1}(\cdot,\cdot)$ is the corresponding energy function (see \S3).
The right hand side of (1.1) is nothing but
the 1DCS, restricted to the highest
weight elements. Hence by the theorem of [JMMO, KMN1] we have
\vskip5mm

\noindent
{\bf Theorem 1.1.} \  For a dominant integral weight
$\la=(\la_1,\cdots,\la_n)$ $(\la_1\geq\cdots\geq\la_n\geq0)$ of $sl_n$,
let us set
$$
V(k\La_0:\lambda)=\{v\in V(k\La_0)\ \vert \ e_jv=0( 1\leq j\leq n-1),
\ \hbox{wt}v=\lambda\}.
$$
We assume $\vert\lambda\vert\equiv 0\hbox{ mod } n$.
Then we have
$$
\hbox{tr}_{V(k\La_0:\lambda)}(q^{-d})
=\lim_{N\ra\infty}q^{-{1\over2}knN(N-1)}K_{\la^{(N)}\mu^{(N)}}(q),
$$
where $\la^{(N)}=
\big(\la_1+kN-{1\over n}\vert\lambda\vert,\cdots,
\la_n+kN-{1\over n}\vert\lambda\vert\big)$
and
$\mu^{(N)}=\big(k^{nN}\big)$.
\vskip5mm

\noindent
This formula was conjectured by A.N. Kirillov (Conjecture 4 in [Kir]).
\vskip3mm

The present paper is organized as follows.
In \S2 we briefly review
the necessary results in crystal theory and
its relation to combinatorics of Young tableaux, the charge, the index, and
the Kostka polynomials.
More precise definitions of crystals and related notation
are given in the Appendix.
We define the energy function and its explicit formula
in \S3. In \S4 the statement of the main theorem of this paper
and its applications
are given. \S5 is devoted to the proof of this main theorem.

\vfill\eject

\noindent
{\fontB 2. Brief Review of Crystal Theory}
\vskip3mm
\noindent
Here we recapitulate some basic facts concerning the crystal base theory
for the case of symmetric and anti-symmetric representations in $sl_n$
mainly through examples.
For more details and precise definitions see the Appendix and [K3,KN,KMN1].

\vskip5mm
\noindent
{\bf 2.1. Crystal base}
\par
\noindent
For a dominant integral weight $\la$ of $sl_n$,
a crystal (base) $B_{\lambda}$ parametrizes a base
of the irreducible representation $V_{\lambda}$.
$B_{\lambda}$ can be thought as
a set of semi-standard tableaux of shape $\lambda$ [KN].
For the symmetric $B_{k \Lambda_1}$
and anti-symmetric $B_{\Lambda_k}$ cases $(1 \leq k \leq n)$,
we also use another representation, as seen in the following
examples. (There are two reasons for this:
(1) We will use semi-standard tableaux to represent a "path".
(2) The energy function has a simple description in this representation.)

\vskip5mm
\noindent
{\bf Example 2.2.}
Elements in $B_{\Lambda_2}=B_{\young(1,1)}$ for $sl_4$ are
$$
	\normalbaselines\mmmth\offinterlineskip
	\vcenter{
   \hbox{$\Flect(0.7cm,0.7cm,\hbox{$\bullet$})$
         }\vskip-0.4pt
   \hbox{$\Flect(0.7cm,0.7cm,\hbox{$\bullet$})$
         }\vskip-0.4pt
   \hbox{$\Flect(0.7cm,0.7cm,\hbox{$$})$
         }\vskip-0.4pt
   \hbox{$\Flect(0.7cm,0.7cm,\hbox{$$})$
         }
                }
(=
	\vcenter{
   \hbox{$\Flect(0.7cm,0.7cm,\hbox{$1$})$
         }\vskip-0.4pt
   \hbox{$\Flect(0.7cm,0.7cm,\hbox{$2$})$
         }
                }),
\hskip1cm
	\vcenter{
   \hbox{$\Flect(0.7cm,0.7cm,\hbox{$\bullet$})$
         }\vskip-0.4pt
   \hbox{$\Flect(0.7cm,0.7cm,\hbox{$$})$
         }\vskip-0.4pt
   \hbox{$\Flect(0.7cm,0.7cm,\hbox{$\bullet$})$
         }\vskip-0.4pt
   \hbox{$\Flect(0.7cm,0.7cm,\hbox{$$})$
         }
                },
\hskip1cm
	\vcenter{
   \hbox{$\Flect(0.7cm,0.7cm,\hbox{$\bullet$})$
         }\vskip-0.4pt
   \hbox{$\Flect(0.7cm,0.7cm,\hbox{$$})$
         }\vskip-0.4pt
   \hbox{$\Flect(0.7cm,0.7cm,\hbox{$$})$
         }\vskip-0.4pt
   \hbox{$\Flect(0.7cm,0.7cm,\hbox{$\bullet$})$
         }
                },
\hskip1cm
	\vcenter{
   \hbox{$\Flect(0.7cm,0.7cm,\hbox{$$})$
         }\vskip-0.4pt
   \hbox{$\Flect(0.7cm,0.7cm,\hbox{$\bullet$})$
         }\vskip-0.4pt
   \hbox{$\Flect(0.7cm,0.7cm,\hbox{$\bullet$})$
         }\vskip-0.4pt
   \hbox{$\Flect(0.7cm,0.7cm,\hbox{$$})$
         }
                },
\hskip1cm
	\vcenter{
   \hbox{$\Flect(0.7cm,0.7cm,\hbox{$$})$
         }\vskip-0.4pt
   \hbox{$\Flect(0.7cm,0.7cm,\hbox{$\bullet$})$
         }\vskip-0.4pt
   \hbox{$\Flect(0.7cm,0.7cm,\hbox{$$})$
         }\vskip-0.4pt
   \hbox{$\Flect(0.7cm,0.7cm,\hbox{$\bullet$})$
         }
                },
\hskip1cm
	\vcenter{
   \hbox{$\Flect(0.7cm,0.7cm,\hbox{$$})$
         }\vskip-0.4pt
   \hbox{$\Flect(0.7cm,0.7cm,\hbox{$$})$
         }\vskip-0.4pt
   \hbox{$\Flect(0.7cm,0.7cm,\hbox{$\bullet$})$
         }\vskip-0.4pt
   \hbox{$\Flect(0.7cm,0.7cm,\hbox{$\bullet$})$
         }
                }.
$$

In general there are $k$-dots in $n$-boxes for $B_{\La_k}$
in the $sl_n$ case.
The boxes are labeled $1,2,3,...,n$ starting from the top,
and each has a weight $\epsilon_i=\Lambda_i-\Lambda_{i-1}$
$(\Lambda_0=\Lambda_{n}=0)$.

\vskip5mm
\noindent
{\bf Example 2.3.}
There are four more elements in $B_{2 \Lambda_1}=B_{\young(2)}$ for $sl_4$ :
$$
	\normalbaselines\mmmth\offinterlineskip
	\vcenter{
   \hbox{$\Flect(0.7cm,0.7cm,\hbox{$\bullet$})$
         }\vskip-0.4pt
   \hbox{$\Flect(0.7cm,0.7cm,\hbox{$\bullet$})$
         }\vskip-0.4pt
   \hbox{$\Flect(0.7cm,0.7cm,\hbox{$$})$
         }\vskip-0.4pt
   \hbox{$\Flect(0.7cm,0.7cm,\hbox{$$})$
         }
                }
(=
	\vcenter{
   \hbox{$\Flect(0.7cm,0.7cm,\hbox{$1$})$\hskip-0.4pt
         $\Flect(0.7cm,0.7cm,\hbox{$2$})$
         }
                }),
\hskip5mm
\hbox{$5$ more and}
\hskip5mm
	\vcenter{
   \hbox{$\Flect(0.7cm,0.7cm,\hbox{$\bullet  \bullet$})$
         }\vskip-0.4pt
   \hbox{$\Flect(0.7cm,0.7cm,\hbox{$$})$
         }\vskip-0.4pt
   \hbox{$\Flect(0.7cm,0.7cm,\hbox{$$})$
         }\vskip-0.4pt
   \hbox{$\Flect(0.7cm,0.7cm,\hbox{$$})$
         }
                },
\hskip1cm
	\vcenter{
   \hbox{$\Flect(0.7cm,0.7cm,\hbox{$$})$
         }\vskip-0.4pt
   \hbox{$\Flect(0.7cm,0.7cm,\hbox{$\bullet \bullet$})$
         }\vskip-0.4pt
   \hbox{$\Flect(0.7cm,0.7cm,\hbox{$$})$
         }\vskip-0.4pt
   \hbox{$\Flect(0.7cm,0.7cm,\hbox{$$})$
         }
                },
\hskip1cm
	\vcenter{
   \hbox{$\Flect(0.7cm,0.7cm,\hbox{$$})$
         }\vskip-0.4pt
   \hbox{$\Flect(0.7cm,0.7cm,\hbox{$$})$
         }\vskip-0.4pt
   \hbox{$\Flect(0.7cm,0.7cm,\hbox{$\bullet \bullet$})$
         }\vskip-0.4pt
   \hbox{$\Flect(0.7cm,0.7cm,\hbox{$$})$
         }
                },
\hskip1cm
	\vcenter{
   \hbox{$\Flect(0.7cm,0.7cm,\hbox{$$})$
         }\vskip-0.4pt
   \hbox{$\Flect(0.7cm,0.7cm,\hbox{$$})$
         }\vskip-0.4pt
   \hbox{$\Flect(0.7cm,0.7cm,\hbox{$$})$
         }\vskip-0.4pt
   \hbox{$\Flect(0.7cm,0.7cm,\hbox{$\bullet \bullet$})$
         }
                }.
$$
\noindent
Note that for $B_{\La_k}$, putting more than one dot in a box is forbidden,
while for $B_{k\La_1}$ it is permitted.
\vskip5mm
\noindent
{\bf 2.4. The action of $\tilde{f_i}$ and $\tilde{e_i}$}
\par
\noindent
On $B_{\lambda}$ one can define the
action of the operators $\tilde{e_i}$ and $\tilde{f_i}$
($ 1 \leq i \leq n-1$) $: B_{\lambda} \rightarrow B_{\lambda} \sqcup \{0\}$.
This is the $q=0$ analog of the action of $sl_n$ generators.
On $B_{\Lambda_k}$ and $B_{k \Lambda_1}$
the action can be described explicitly in the following manners [KN].

{\bf Case 1: \  $b \in B_{\Lambda_k}$}.

(a) If the $i$-th box in $b$ is dotted and the $(i+1)$-th is not,
then $\tilde{f_i}b$
is the element obtained by shifting the $i$-th dot to $(i+1)$-th box.
Otherwise, $\tilde{f_i}b=0$.

(b) If the $(i+1)$-th box in $b$ is dotted and the $i$-th is not,
then $\tilde{e_i}b$
is the element obtained by shifting the $(i+1)$-th dot to $i$-th box.
Otherwise, $\tilde{e_i}b=0$.

{\bf Case 2: \  $b \in B_{k \Lambda_1}$}.

(a) If there are dots in the $i$-th box in $b$, then $\tilde{f_i}b$
is the element obtained by shifting one of the $i$-th dots to
the $(i+1)$-th box.
Otherwise, $\tilde{f_i}b=0$.

(b) If there are dots in the $(i+1)$-th box in $b$, then $\tilde{e_i}b$
is the element obtained by shifting one of the $(i+1)$-th dots to the
$i$-th box.
Otherwise, $\tilde{e_i}b=0$.

\vskip5mm
One can also define the action of $\tilde{e_0}$ and $\tilde{f_0}$,
which we identify with $\tilde{e_n}$ and $\tilde{f_n}$.
Here the column is considered to be periodically continued with period
$n$. Then the action of $\tilde{e_n},\tilde{f_n}$ is defined by the above
rules.
In this way $B_{\lambda}$ is extended to a crystal for
$\widehat{sl_n}$ [KMN2]. This extension is crucial for
the characterization of the energy function.

\vskip5mm
\noindent
{\bf 2.5. Tensor products}
\par
\noindent
On the tensor product of crystals (direct product as a set),
the actions of $\tilde{e_i}$ and $\tilde{f_i}$ are defined as follows.
Since the action depends only on the configuration in the $i$-th and
$(i+1)$-th box, only these two boxes will be considered.

In both symmetric and anti-symmetric cases, the rules are
as follows:

(a) We neglect empty boxes $\Young(1,1)$.

(b) We remove singlet pairs of dots, that is
$$
	\normalbaselines\mmmth\offinterlineskip
	\vcenter{
   \hbox{$\Flect(0.7cm,0.7cm,\hbox{$(\bullet)$})$
         }\vskip-0.4pt
   \hbox{$\Flect(0.7cm,0.7cm,\hbox{$(\bullet)$})$
         }
                }
\ {\rm or} \
	\vcenter{
   \hbox{$\Flect(0.7cm,0.7cm,\hbox{$(\bullet)$})$
         }\vskip-0.4pt
   \hbox{$\Flect(0.7cm,0.7cm,\hbox{$$})$
         }
                }
\otimes
	\vcenter{
   \hbox{$\Flect(0.7cm,0.7cm,\hbox{$$})$
         }\vskip-0.4pt
   \hbox{$\Flect(0.7cm,0.7cm,\hbox{$(\bullet)$})$
         }
                }
$$
for the anti-symmetric case and
$$
	\normalbaselines\mmmth\offinterlineskip
	\vcenter{
   \hbox{$\Flect(0.7cm,1.5cm,\hbox{$\bullet (\bullet \bullet \bullet)$})$
         }\vskip-0.4pt
   \hbox{$\Flect(0.7cm,1.5cm,\hbox{any})$
         }
                }
\otimes
	\vcenter{
   \hbox{$\Flect(0.7cm,1.5cm,\hbox{any})$
         }\vskip-0.4pt
   \hbox{$\Flect(0.7cm,1.5cm,\hbox{$(\bullet \bullet \bullet)$})$
         }
                }
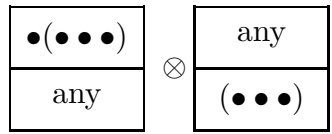
$$
(and others like this) for the symmetric case.

(c) We continue (a) and (b) as far as possible until we obtain
something like
$$
	\normalbaselines\mmmth\offinterlineskip
	\vcenter{
   \hbox{$\Flect(0.7cm,0.7cm,\hbox{$$})$
         }\vskip-0.4pt
   \hbox{$\Flect(0.7cm,0.7cm,\hbox{$\bullet$})$
         }
                }
\otimes
	\vcenter{
   \hbox{$\Flect(0.7cm,0.7cm,\hbox{$\bullet$})$
         }\vskip-0.4pt
   \hbox{$\Flect(0.7cm,0.7cm,\hbox{$$})$
         }
                }
\otimes
	\vcenter{
   \hbox{$\Flect(0.7cm,0.7cm,\hbox{$\bullet$})$
         }\vskip-0.4pt
   \hbox{$\Flect(0.7cm,0.7cm,\hbox{$$})$
         }
                }
$$
for the anti-symmetric case and
$$
	\normalbaselines\mmmth\offinterlineskip
	\vcenter{
   \hbox{$\Flect(0.7cm,0.7cm,\hbox{$\bullet$})$
         }\vskip-0.4pt
   \hbox{$\Flect(0.7cm,0.7cm,\hbox{$\bullet \bullet$})$
         }
                }
\otimes
	\vcenter{
   \hbox{$\Flect(0.7cm,0.7cm,\hbox{$\bullet$})$
         }\vskip-0.4pt
   \hbox{$\Flect(0.7cm,0.7cm,\hbox{$$})$
         }
                }
\otimes
	\vcenter{
   \hbox{$\Flect(0.7cm,0.7cm,\hbox{$\bullet \bullet$})$
         }\vskip-0.4pt
   \hbox{$\Flect(0.7cm,0.7cm,\hbox{$$})$
         }
                }
$$
for the symmetric case.

(d) Then $\tilde{e_i}$ shifts a dot in the bottom box of the rightmost
to the top box
and $\tilde{f_i}$ shifts a dot in the top box of the leftmost
to the bottom box.
If there is no such dot, then $\tilde{e_i} b=0$ (resp. $\tilde{f_i}b=0$).

\vskip5mm
\noindent
{\bf Example 2.6.}
Let $b \in B_{\Lambda_3} \otimes B_{\Lambda_2} \otimes B_{\Lambda_2}
\otimes B_{\Lambda_1} \otimes B_{\Lambda_1}\otimes B_{\Lambda_1}$ as
follows :
$$
	\normalbaselines\mmmth\offinterlineskip
	\vcenter{
   \hbox{$\Flect(0.7cm,0.7cm,\hbox{$\bullet$})$
         }\vskip-0.4pt
   \hbox{$\Flect(0.7cm,0.7cm,\hbox{$\bullet$})$
         }\vskip-0.4pt
   \hbox{$\Flect(0.7cm,0.7cm,\hbox{$\bullet$})$
         }\vskip-0.4pt
   \hbox{$\Flect(0.7cm,0.7cm,\hbox{$$})$
         }
                }
\otimes
	\vcenter{
   \hbox{$\Flect(0.7cm,0.7cm,\hbox{$\bullet$})$
         }\vskip-0.4pt
   \hbox{$\Flect(0.7cm,0.7cm,\hbox{$$})$
         }\vskip-0.4pt
   \hbox{$\Flect(0.7cm,0.7cm,\hbox{$$})$
         }\vskip-0.4pt
   \hbox{$\Flect(0.7cm,0.7cm,\hbox{$\bullet$})$
         }
                }
\otimes
	\vcenter{
   \hbox{$\Flect(0.7cm,0.7cm,\hbox{$\bullet$})$
         }\vskip-0.4pt
   \hbox{$\Flect(0.7cm,0.7cm,\hbox{$$})$
         }\vskip-0.4pt
   \hbox{$\Flect(0.7cm,0.7cm,\hbox{$\bullet$})$
         }\vskip-0.4pt
   \hbox{$\Flect(0.7cm,0.7cm,\hbox{$$})$
         }
                }
\otimes
	\vcenter{
   \hbox{$\Flect(0.7cm,0.7cm,\hbox{$$})$
         }\vskip-0.4pt
   \hbox{$\Flect(0.7cm,0.7cm,\hbox{$\bullet$})$
         }\vskip-0.4pt
   \hbox{$\Flect(0.7cm,0.7cm,\hbox{$$})$
         }\vskip-0.4pt
   \hbox{$\Flect(0.7cm,0.7cm,\hbox{$$})$
         }
                }
\otimes
	\vcenter{
   \hbox{$\Flect(0.7cm,0.7cm,\hbox{$$})$
         }\vskip-0.4pt
   \hbox{$\Flect(0.7cm,0.7cm,\hbox{$\bullet$})$
         }\vskip-0.4pt
   \hbox{$\Flect(0.7cm,0.7cm,\hbox{$$})$
         }\vskip-0.4pt
   \hbox{$\Flect(0.7cm,0.7cm,\hbox{$$})$
         }
                }
\otimes
	\vcenter{
   \hbox{$\Flect(0.7cm,0.7cm,\hbox{$$})$
         }\vskip-0.4pt
   \hbox{$\Flect(0.7cm,0.7cm,\hbox{$\bullet$})$
         }\vskip-0.4pt
   \hbox{$\Flect(0.7cm,0.7cm,\hbox{$$})$
         }\vskip-0.4pt
   \hbox{$\Flect(0.7cm,0.7cm,\hbox{$$})$
         }
                }.
$$
As far as the actions of $\tilde{e_1}$ and $\tilde{f_1}$ are concerned,
$b$ is equivalent to
$$
	\normalbaselines\mmmth\offinterlineskip
	\vcenter{
   \hbox{$\Flect(0.7cm,0.7cm,\hbox{$\bullet$})$
         }\vskip-0.4pt
   \hbox{$\Flect(0.7cm,0.7cm,\hbox{$\bullet$})$
         }
                }
\otimes
	\vcenter{
   \hbox{$\Flect(0.7cm,0.7cm,\hbox{$\bullet$})$
         }\vskip-0.4pt
   \hbox{$\Flect(0.7cm,0.7cm,\hbox{$$})$
         }
                }
\otimes
	\vcenter{
   \hbox{$\Flect(0.7cm,0.7cm,\hbox{$\bullet$})$
         }\vskip-0.4pt
   \hbox{$\Flect(0.7cm,0.7cm,\hbox{$$})$
         }
                }
\otimes
	\vcenter{
   \hbox{$\Flect(0.7cm,0.7cm,\hbox{$$})$
         }\vskip-0.4pt
   \hbox{$\Flect(0.7cm,0.7cm,\hbox{$\bullet$})$
         }
                }
\otimes
	\vcenter{
   \hbox{$\Flect(0.7cm,0.7cm,\hbox{$$})$
         }\vskip-0.4pt
   \hbox{$\Flect(0.7cm,0.7cm,\hbox{$\bullet$})$
         }
                }
\otimes
	\vcenter{
   \hbox{$\Flect(0.7cm,0.7cm,\hbox{$$})$
         }\vskip-0.4pt
   \hbox{$\Flect(0.7cm,0.7cm,\hbox{$\bullet$})$
         }
                }.
$$
Then by the above rules (a) and (b), this can be reduced to the following:
$$
	\normalbaselines\mmmth\offinterlineskip
\rightarrow
	\vcenter{
   \hbox{$\Flect(0.7cm,0.7cm,\hbox{$$})$
         }\vskip-0.4pt
   \hbox{$\Flect(0.7cm,0.7cm,\hbox{$$})$
         }
                }
\otimes
	\vcenter{
   \hbox{$\Flect(0.7cm,0.7cm,\hbox{$\bullet$})$
         }\vskip-0.4pt
   \hbox{$\Flect(0.7cm,0.7cm,\hbox{$$})$
         }
                }
\otimes
	\vcenter{
   \hbox{$\Flect(0.7cm,0.7cm,\hbox{$$})$
         }\vskip-0.4pt
   \hbox{$\Flect(0.7cm,0.7cm,\hbox{$$})$
         }
                }
\otimes
	\vcenter{
   \hbox{$\Flect(0.7cm,0.7cm,\hbox{$$})$
         }\vskip-0.4pt
   \hbox{$\Flect(0.7cm,0.7cm,\hbox{$$})$
         }
                }
\otimes
	\vcenter{
   \hbox{$\Flect(0.7cm,0.7cm,\hbox{$$})$
         }\vskip-0.4pt
   \hbox{$\Flect(0.7cm,0.7cm,\hbox{$\bullet$})$
         }
                }
\otimes
	\vcenter{
   \hbox{$\Flect(0.7cm,0.7cm,\hbox{$$})$
         }\vskip-0.4pt
   \hbox{$\Flect(0.7cm,0.7cm,\hbox{$\bullet$})$
         }
                }
\hskip1cm
{\rm by \ (b)}
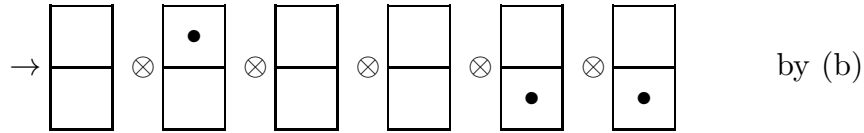
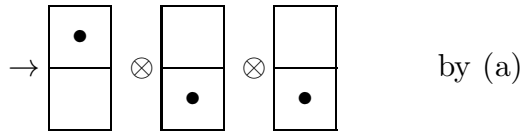
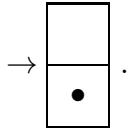
$$
and
$$
	\normalbaselines\mmmth\offinterlineskip
\rightarrow
	\vcenter{
   \hbox{$\Flect(0.7cm,0.7cm,\hbox{$\bullet$})$
         }\vskip-0.4pt
   \hbox{$\Flect(0.7cm,0.7cm,\hbox{$$})$
         }
                }
\otimes
	\vcenter{
   \hbox{$\Flect(0.7cm,0.7cm,\hbox{$$})$
         }\vskip-0.4pt
   \hbox{$\Flect(0.7cm,0.7cm,\hbox{$\bullet$})$
         }
                }
\otimes
	\vcenter{
   \hbox{$\Flect(0.7cm,0.7cm,\hbox{$$})$
         }\vskip-0.4pt
   \hbox{$\Flect(0.7cm,0.7cm,\hbox{$\bullet$})$
         }
                }
\hskip1cm
{\rm by \ (a)}
$$
and finally
$$
	\normalbaselines\mmmth\offinterlineskip
\rightarrow
	\vcenter{
   \hbox{$\Flect(0.7cm,0.7cm,\hbox{$$})$
         }\vskip-0.4pt
   \hbox{$\Flect(0.7cm,0.7cm,\hbox{$\bullet$})$
         }
                }.
$$
Hence $\tilde{e_1}$ acts on the last component of the tensor product,
while $\tilde{f_1}b=0$.

\vskip5mm
\noindent
{\bf 2.7. Parametrization of the highest weight vectors}
\par
\noindent
In what follows we fix a partition $\mu=(\mu_1,\cdots,\mu_l)$
and consider the tensor products
$$
{\cal H}_{\mu}=
B_{\mu_1 \Lambda_1} \otimes \cdots \otimes B_{\mu_l \Lambda_1},
$$
and
$$
{\cal H'}_{\mu}=
B_{\Lambda_{\mu_1}} \otimes \cdots \otimes B_{\Lambda_{\mu_l}}.
$$

Here we will give a canonical parametrization of the $sl_n$ highest
weight elements $b$ (i.e. $\tilde{e_i}b=0$ for all $1\leq i\leq n-1$)
in terms of semi-standard tableau.

\vskip5mm
\noindent
{\bf 2.8. Semi-standard tableau ${\cal T}(\lambda,\mu)$}
\par
\noindent
For partitions $\lambda$ and $\mu$ ($\vert \lambda \vert = \vert \mu \vert$),
let ${\cal T}(\lambda, \mu)$ be the set of semi-standard tableaux
of shape $\lambda$ and weight $\mu$.

\vskip5mm
\noindent
{\bf Example 2.9.} \
A semi-standard tableau of shape $\lambda=(5,4,2)$ is
a diagram of the form
$$
	\normalbaselines\mmmth\offinterlineskip
	\vcenter{
   \hbox{$\Flect(0.7cm,0.7cm,\hbox{$i_{1,1}$})\hskip-0.4pt
          \Flect(0.7cm,0.7cm,\hbox{$i_{1,2}$})\hskip-0.4pt
          \Flect(0.7cm,0.7cm,\hbox{$i_{1,3}$})\hskip-0.4pt
          \Flect(0.7cm,0.7cm,\hbox{$i_{1,4}$})\hskip-0.4pt
          \Flect(0.7cm,0.7cm,\hbox{$i_{1,5}$})$
         }\vskip-0.4pt
   \hbox{$\Flect(0.7cm,0.7cm,\hbox{$i_{2,1}$})\hskip-0.4pt
          \Flect(0.7cm,0.7cm,\hbox{$i_{2,2}$})\hskip-0.4pt
          \Flect(0.7cm,0.7cm,\hbox{$i_{2,3}$})\hskip-0.4pt
          \Flect(0.7cm,0.7cm,\hbox{$i_{2,4}$})$
         }\vskip-0.4pt
   \hbox{$\Flect(0.7cm,0.7cm,\hbox{$i_{3,1}$})\hskip-0.4pt
          \Flect(0.7cm,0.7cm,\hbox{$i_{3,2}$})$
         }
                }
$$
where the numbers $i_{n,m}$ satisfy the following conditions:
$$
i_{n,m+1} \geq i_{n,m},
\hskip1cm
i_{n+1,m} > i_{n,m}.
$$

\vskip5mm
\noindent
{\bf Example 2.10.} \
The following tableau $T$ has the shape $\lambda=(5,4,2)$ and
the weight $\mu=(3^2,2^2,1)$.
$$
	\normalbaselines\mmmth\offinterlineskip
	\vcenter{
   \hbox{$\Flect(0.7cm,0.7cm,\hbox{$1$})\hskip-0.4pt
          \Flect(0.7cm,0.7cm,\hbox{$1$})\hskip-0.4pt
          \Flect(0.7cm,0.7cm,\hbox{$1$})\hskip-0.4pt
          \Flect(0.7cm,0.7cm,\hbox{$2$})\hskip-0.4pt
          \Flect(0.7cm,0.7cm,\hbox{$4$})$
         }\vskip-0.4pt
   \hbox{$\Flect(0.7cm,0.7cm,\hbox{$2$})\hskip-0.4pt
          \Flect(0.7cm,0.7cm,\hbox{$2$})\hskip-0.4pt
          \Flect(0.7cm,0.7cm,\hbox{$3$})\hskip-0.4pt
          \Flect(0.7cm,0.7cm,\hbox{$4$})$
         }\vskip-0.4pt
   \hbox{$\Flect(0.7cm,0.7cm,\hbox{$3$})\hskip-0.4pt
          \Flect(0.7cm,0.7cm,\hbox{$5$})$
         }
                }
$$
Here, $\sharp \{ (n,m) \ \vert \ i_{n,m}=j \}=\mu_j$.
\vskip5mm

\noindent
We denote by $\la^\prime$ the transposition of the partition $\la$ [M].
For each $T \in {\cal T}(\lambda,\mu)$
one can construct a highest weight element $b$ in ${\cal H}_{\mu}$
[resp. ${\cal H'}_{\mu}$] of weight $\lambda$ [resp. $\lambda'$]
as in the following examples.

\vskip5mm
\noindent
{\bf Example 2.11.} \

\noindent
The highest weight element in ${\cal H}_{\mu}$
corresponding to the tableau $T$ in
Example 2.10 is
$$
	\normalbaselines\mmmth\offinterlineskip
	\vcenter{
   \hbox{$\Flect(0.7cm,0.9cm,\hbox{$\bullet \bullet \bullet$})$
         }\vskip-0.4pt
   \hbox{$\Flect(0.7cm,0.9cm,\hbox{$$})$
         }\vskip-0.4pt
   \hbox{$\Flect(0.7cm,0.9cm,\hbox{$$})$
         }\vskip-0.4pt
   \hbox{$\Flect(0.7cm,0.9cm,\hbox{$$})$
         }
                }
\otimes
	\vcenter{
   \hbox{$\Flect(0.7cm,0.7cm,\hbox{$\bullet$})$
         }\vskip-0.4pt
   \hbox{$\Flect(0.7cm,0.7cm,\hbox{$\bullet \bullet$})$
         }\vskip-0.4pt
   \hbox{$\Flect(0.7cm,0.7cm,\hbox{$$})$
         }\vskip-0.4pt
   \hbox{$\Flect(0.7cm,0.7cm,\hbox{$$})$
         }
                }
\otimes
	\vcenter{
   \hbox{$\Flect(0.7cm,0.7cm,\hbox{$$})$
         }\vskip-0.4pt
   \hbox{$\Flect(0.7cm,0.7cm,\hbox{$\bullet$})$
         }\vskip-0.4pt
   \hbox{$\Flect(0.7cm,0.7cm,\hbox{$\bullet$})$
         }\vskip-0.4pt
   \hbox{$\Flect(0.7cm,0.7cm,\hbox{$$})$
         }
                }
\otimes
	\vcenter{
   \hbox{$\Flect(0.7cm,0.7cm,\hbox{$\bullet$})$
         }\vskip-0.4pt
   \hbox{$\Flect(0.7cm,0.7cm,\hbox{$\bullet$})$
         }\vskip-0.4pt
   \hbox{$\Flect(0.7cm,0.7cm,\hbox{$$})$
         }\vskip-0.4pt
   \hbox{$\Flect(0.7cm,0.7cm,\hbox{$$})$
         }
                }
\otimes
	\vcenter{
   \hbox{$\Flect(0.7cm,0.7cm,\hbox{$$})$
         }\vskip-0.4pt
   \hbox{$\Flect(0.7cm,0.7cm,\hbox{$$})$
         }\vskip-0.4pt
   \hbox{$\Flect(0.7cm,0.7cm,\hbox{$\bullet$})$
         }\vskip-0.4pt
   \hbox{$\Flect(0.7cm,0.7cm,\hbox{$$})$
         }

                }.
$$

\noindent
The highest weight element in ${\cal H'}_{\mu}$
corresponding to the tableau $T$ in
Example 2.10 is
$$
	\normalbaselines\mmmth\offinterlineskip
	\vcenter{
   \hbox{$\Flect(0.7cm,0.7cm,\hbox{$\bullet$})$
         }\vskip-0.4pt
   \hbox{$\Flect(0.7cm,0.7cm,\hbox{$\bullet$})$
         }\vskip-0.4pt
   \hbox{$\Flect(0.7cm,0.7cm,\hbox{$\bullet$})$
         }\vskip-0.4pt
   \hbox{$\Flect(0.7cm,0.7cm,\hbox{$$})$
         }\vskip-0.4pt
   \hbox{$\Flect(0.7cm,0.7cm,\hbox{$$})$
         }
                }
\otimes
	\vcenter{
   \hbox{$\Flect(0.7cm,0.7cm,\hbox{$\bullet$})$
         }\vskip-0.4pt
   \hbox{$\Flect(0.7cm,0.7cm,\hbox{$\bullet$})$
         }\vskip-0.4pt
   \hbox{$\Flect(0.7cm,0.7cm,\hbox{$$})$
         }\vskip-0.4pt
   \hbox{$\Flect(0.7cm,0.7cm,\hbox{$\bullet$})$
         }\vskip-0.4pt
   \hbox{$\Flect(0.7cm,0.7cm,\hbox{$$})$
         }
                }
\otimes
	\vcenter{
   \hbox{$\Flect(0.7cm,0.7cm,\hbox{$\bullet$})$
         }\vskip-0.4pt
   \hbox{$\Flect(0.7cm,0.7cm,\hbox{$$})$
         }\vskip-0.4pt
   \hbox{$\Flect(0.7cm,0.7cm,\hbox{$\bullet$})$
         }\vskip-0.4pt
   \hbox{$\Flect(0.7cm,0.7cm,\hbox{$$})$
         }\vskip-0.4pt
   \hbox{$\Flect(0.7cm,0.7cm,\hbox{$$})$
         }
                }
\otimes
	\vcenter{
   \hbox{$\Flect(0.7cm,0.7cm,\hbox{$$})$
         }\vskip-0.4pt
   \hbox{$\Flect(0.7cm,0.7cm,\hbox{$$})$
         }\vskip-0.4pt
   \hbox{$\Flect(0.7cm,0.7cm,\hbox{$$})$
         }\vskip-0.4pt
   \hbox{$\Flect(0.7cm,0.7cm,\hbox{$\bullet$})$
         }\vskip-0.4pt
   \hbox{$\Flect(0.7cm,0.7cm,\hbox{$\bullet$})$
         }
                }
\otimes
	\vcenter{
   \hbox{$\Flect(0.7cm,0.7cm,\hbox{$$})$
         }\vskip-0.4pt
   \hbox{$\Flect(0.7cm,0.7cm,\hbox{$\bullet$})$
         }\vskip-0.4pt
   \hbox{$\Flect(0.7cm,0.7cm,\hbox{$$})$
         }\vskip-0.4pt
   \hbox{$\Flect(0.7cm,0.7cm,\hbox{$$})$
         }\vskip-0.4pt
   \hbox{$\Flect(0.7cm,0.7cm,\hbox{$$})$
         }

                }.
$$
Note that the correspondence makes sense iff
$n \geq l(\lambda)$ [resp. $n \geq l(\lambda')$].

\vskip5mm
\noindent
{\bf Theorem 2.12.} \
The element $b \in {\cal H}_{\mu}$ [resp. ${\cal H'}_{\mu}$]
is a highest weight element with weight $\lambda$ [resp. $\lambda'$]
if and only if
it is an image of $T \in {\cal T}(\lambda,\mu)$ under the
above correspondence.
This correspondence gives a bijection between the highest elements
and the semi-standard tableaux.
\noindent
Proof. The statement is easily proved using the rules
of the action of $\tilde{e_i}$ described above. q.e.d.

\vskip1cm
\noindent
{\bf 2.13. Review of the Kostka Polynomials}
\vskip2mm
\noindent
Here we review the necessary basic facts concerning the symmetric functions,
in particular Kostka polynomials. For more
details, see Ref.[M].

\vskip3mm
\noindent
{\bf 2.14. Kostka number $K_{\lambda,\mu}$}
\par
\noindent
For a partition $\lambda$,
let $s_{\lambda}(x)$, $m_{\lambda}(x)$, $h_{\lambda}(x)$
and $e_{\lambda}(x)$ be the symmetric functions on $n$-variables
$x=(x_1,\cdots,x_n)$.
These are referred as the Schur, monomial, complete and elementary
functions respectively. We are interested in their transition matrices.

The Kostka number $K_{\lambda,\mu}$ is defined by
$$
s_{\lambda}(x)=\sum_{\mu} K_{\lambda,\mu} m_{\mu}(x),
$$
that is, $K_{\lambda,\mu}$ is the multiplicity of weight $\mu$
in the irreducible highest weight representation $V_{\lambda}$.

By the duality relation
$$
\sum_{\lambda} m_{\lambda}(x) h_{\lambda}(y)=
\sum_{\lambda} s_{\lambda}(x) s_{\lambda}(y),
$$
this definition is equivalent to
$$
h_{\mu}(x)=\sum_{\lambda} K_{\lambda, \mu} s_{\lambda}(x).
$$
Since $h_{\mu}=h_{\mu_1} \cdots h_{\mu_l}$ is the character of
the tensor product of complete symmetric representations
$V_{\mu_1 \Lambda_1} \otimes \cdots \otimes V_{\mu_l \Lambda_1}$,
$K_{\lambda, \mu}$ is the multiplicity of $V_{\lambda}$ in this
tensor product, i.e., the Clebsch-Gordan coefficient.

By applying the involution $\omega$, $\omega(s_{\lambda})=s_{\lambda'}$,
for which $\omega(h_\mu)=e_\mu$ to the above equation, we also have
$$
e_{\mu}(x)=\sum_{\lambda} K_{\lambda,\mu} s_{\lambda'}(x),
$$
that is, $K_{\lambda, \mu}$ is the multiplicity of $V_{\lambda'}$
in the tensor products
$V_{\Lambda_{\mu_1}} \otimes \cdots \otimes V_{\Lambda_{\mu_l}}$.

Theorem 2.12 gives an
explanation in terms of crystals
for the well known fact that $K_{\lambda,\mu}$ is the number of
semi-standard tableaux of shape $\lambda$ and weight $\mu$.

\vskip1cm
\noindent
{\bf 2.15. Kostka polynomial $K_{\lambda,\mu}(q)$}
\par
\noindent
The Kostka polynomial $K_{\lambda, \mu}(q)$ is defined
as the transition coefficient of the Schur function
in terms of the Hall-Littlewood polynomial $P_\mu(x,q)$ :
$$
s_{\lambda}(x)=\sum_{\mu} K_{\lambda,\mu}(q) P_{\mu}(x,q).
$$
Since $P_{\lambda}(x,1)=m_{\lambda}(x)$, $K_{\lambda,\mu}(q)$
is the $q$-analog of the Kostka number $K_{\lambda,\mu}=K_{\lambda,\mu}(1)$.

Though it is quite nontrivial from this definition, $K_{\lambda,\mu}(q)$
is a polynomial in $q$ with {\it nonnegative integer} coefficients.
This fact was proved by Lascoux-Schutzenberger by establishing
the following combinatorial formula for Kostka polynomials.

\vskip5mm
\noindent
{\bf Theorem 2.16.} ([LS]) \
$$
K_{\lambda,\mu}(q)=\sum_{T \in {\cal T}(\lambda,\mu)} q^{c(T)},
$$
where the definition of the charge
$c \ : \ {\cal T} \rightarrow {\bf Z}_{\geq 0}$
is illustrated by the following example [LS,M].

\vskip5mm
\noindent
{\bf Example 2.17.} \
For the tableau $T$ in Example 2.10, the charge is calculated
by the following three steps.
$$
	\normalbaselines\mmmth\offinterlineskip
	\vcenter{
   \hbox{$\Flect(0.7cm,0.7cm,\hbox{$1$})\hskip-0.4pt
          \Flect(0.7cm,0.7cm,\hbox{$1$})\hskip-0.4pt
          \Flect(0.7cm,0.7cm,\hbox{${\bf 1}_0$})\hskip-0.4pt
          \Flect(0.7cm,0.7cm,\hbox{$2$})\hskip-0.4pt
          \Flect(0.7cm,0.7cm,\hbox{${\bf 4}_1$})$
         }\vskip-0.4pt
   \hbox{$\Flect(0.7cm,0.7cm,\hbox{$2$})\hskip-0.4pt
          \Flect(0.7cm,0.7cm,\hbox{${\bf 2}_0$})\hskip-0.4pt
          \Flect(0.7cm,0.7cm,\hbox{$3$})\hskip-0.4pt
          \Flect(0.7cm,0.7cm,\hbox{$4$})$
         }\vskip-0.4pt
   \hbox{$\Flect(0.7cm,0.7cm,\hbox{${\bf 3}_0$})\hskip-0.4pt
          \Flect(0.7cm,0.7cm,\hbox{${\bf 5}_1$})$
         }
                }
\rightarrow
	\vcenter{
   \hbox{$\Flect(0.7cm,0.7cm,\hbox{$1$})\hskip-0.4pt
          \Flect(0.7cm,0.7cm,\hbox{${\bf 1}_0$})\hskip-0.4pt
          \Flect(0.7cm,0.7cm,\hbox{$$})\hskip-0.4pt
          \Flect(0.7cm,0.7cm,\hbox{$2$})\hskip-0.4pt
          \Flect(0.7cm,0.7cm,\hbox{$$})$
         }\vskip-0.4pt
   \hbox{$\Flect(0.7cm,0.7cm,\hbox{${\bf 2}_0$})\hskip-0.4pt
          \Flect(0.7cm,0.7cm,\hbox{$$})\hskip-0.4pt
          \Flect(0.7cm,0.7cm,\hbox{${\bf 3}_1$})\hskip-0.4pt
          \Flect(0.7cm,0.7cm,\hbox{${\bf 4}_2$})$
         }\vskip-0.4pt
   \hbox{$\Flect(0.7cm,0.7cm,\hbox{$$})\hskip-0.4pt
          \Flect(0.7cm,0.7cm,\hbox{$$})$
         }
                }
\rightarrow
	\vcenter{
   \hbox{$\Flect(0.7cm,0.7cm,\hbox{${\bf 1}_0$})\hskip-0.4pt
          \Flect(0.7cm,0.7cm,\hbox{$$})\hskip-0.4pt
          \Flect(0.7cm,0.7cm,\hbox{$$})\hskip-0.4pt
          \Flect(0.7cm,0.7cm,\hbox{${\bf 2}_1$})\hskip-0.4pt
          \Flect(0.7cm,0.7cm,\hbox{$$})$
         }\vskip-0.4pt
   \hbox{$\Flect(0.7cm,0.7cm,\hbox{$$})\hskip-0.4pt
          \Flect(0.7cm,0.7cm,\hbox{$$})\hskip-0.4pt
          \Flect(0.7cm,0.7cm,\hbox{$$})\hskip-0.4pt
          \Flect(0.7cm,0.7cm,\hbox{$$})$
         }\vskip-0.4pt
   \hbox{$\Flect(0.7cm,0.7cm,\hbox{$$})\hskip-0.4pt
          \Flect(0.7cm,0.7cm,\hbox{$$})$
         }
                }
$$
In each step, the sequence ${\bf 1},{\bf 2},{\bf 3},...$
is selected from the top and right (periodically),
and each suffix counts the "winding number" (=
the number of times the top is revisited).
Then index, ${\rm ind}(j)$, is defined as a sum of the
winding numbers assigned to ${\bf j}$, i.e.,
$$
{\rm ind}(1)=0, \hskip5mm
{\rm ind}(2)=1, \hskip5mm
{\rm ind}(3)=1, \hskip5mm
{\rm ind}(4)=3, \hskip5mm
{\rm ind}(5)=1.
$$
Finally, the charge is the sum of all indices, $c(T)=6$.

\vskip5mm
In the main Theorem, we will give a similar (but not the same)
formula for the index (or $K_{\lambda,\mu}(q)$) from the point of view
of crystal theory.

\vskip1cm
\noindent
{\bf 2.18. Skew Kostka polynomials}
\par
\noindent
For the skew Young diagram $\lambda / \nu$,
the skew Kostka polynomial
$K_{\lambda / \nu, \mu}(q)$
was similarly defined by [KR] as
the connection matrix between the skew Schur function
$s_{\lambda / \nu}$ and the Hall polynomials.
It has the description [KR, B]
$$
K_{\lambda / \nu, \mu}(q)=\sum_{T \in {\cal T}(\lambda / \nu, \mu)} q^{c(T)},
$$
where ${\cal T}(\lambda / \nu, \mu)$ is the set of semi-standard tableaux
of shape $\lambda / \nu$ and weight $\mu$.

Similarly to the non-skew case ($\nu=\phi$), there exists a one to one
correspondence between the semi-standard tableau $T$ in
${\cal T}(\lambda / \nu, \mu)$ and the highest element $b$ in
$B_{\nu} \otimes {\cal H}_{\mu}$ of weight $\lambda$,
[resp. $b'$ in
$B_{\nu'} \otimes {\cal H}'_{\mu}$ of weight $\lambda'$].
In terms of paths, these correspond to the fusion path
beginning at $\nu$ and
ending at $\lambda$ [resp.  $\nu'$ and $\lambda'$ ].
\vfill\eject

%%%%%%%%%%%%%%%%%%%%%%%%%%%%%%%%%%%%%%%%

%Nakayashiki personal macros
\def\ot{\otimes}
\def\ra{\longrightarrow}
\def\et{\tilde{e}}
\def\ft{\tilde{f}}
\def\La{\Lambda}
\def\la{\lambda}
\def\mapright#1{\smash{\mathop{\longrightarrow}\limits^{#1}}}

\noindent
{\fontB 3. Energy Functions}\par
\vskip5mm
\noindent
{\bf 3.1. General definition}\par
\noindent
Here we give definitions and descriptions
of energy functions.
The precise definitions of notation used here
are explained briefly in the Appendix.

Let us
consider two classical crystals, $B_1$ and $B_2$, for $\widehat{sl_n}$
satisfying the following conditions:
$$
\eqalignno{
&\hbox{$B_1\ot B_2$ is connected.}&(3.1)\cr
&\hbox{$B_1\ot B_2$ is isomorphic to $B_2\ot B_1$.}&(3.2)\cr
}
$$

\vskip5mm
\noindent
{\bf Definition 3.2.}
Suppose that $b_1\ot b_2$ is mapped to $b_2^\prime\ot b_1^\prime$
by the isomorphism (3.2):
$$
\eqalign{
B_1\ot B_2&\simeq B_2\ot B_1,\cr
b_1\ot b_2&\mapsto b_2^\prime\ot b_1^\prime.\cr
}
$$
The energy function $H$ of the crystal $B_1\ot B_2$
is a map
$$
H: B_1\ot B_2\ra {\bf Z}
$$
which satisfies the following conditions.
Suppose that $\et_i(b_1\ot b_2)\neq0$.
Then
$$
\eqalign{
H\big(\et_i(b_1\ot b_2)\big)
&=H(b_1\ot b_2)+1\hbox{ if $i=0$, $\varphi_0(b_1)\geq\varepsilon_0(b_2)$
and $\varphi_0(b_2^\prime)\geq\varepsilon_0(b_1^\prime)$},
\cr
&=H(b_1\ot b_2)-1\hbox{ if $i=0$, $\varphi_0(b_1)<\varepsilon_0(b_2)$
and $\varphi_0(b_2^\prime)<\varepsilon_0(b_1^\prime)$},
\cr
&=H(b_1\ot b_2)\hbox{ otherwise}.
\cr}
$$
We often denote $H(b_1,b_2)$ instead of $H(b_1 \otimes b_2)$
for the sake of simplicity.
\vskip5mm

By definition the energy function is constant on each connected
component of $B_1\ot B_2$ as $sl_n$ crystal
and is unique up to an additive constant by the condition (3.1).
For $B_1=B_2$ the energy function has been used extensively
to evaluate the one point function of solvable lattice
models [DJKMO,KMN1,DJO].
The study of spinon character formulas [NY1-2] leads us to
consider the above slightly generalized energy functions.

\vskip5mm
\noindent
{\bf Remark 3.3.}
(1). The definition of the energy function above follows from the
property of the combinatorial $R$ matrix in a manner similar
to the $B_1=B_2$ case [KMN1].\par
\noindent (2). In [KMN2] the crystal $B_{k\La_l}$
(for $\widehat{sl_n}$) is proved to be a perfect crystal of level
$k$ (Definition 4.6.1 [KMN1]). Thus
if $B_1$ and $B_2$ are crystals of the form
$B_{k\La_l}$, we can prove that
the condition (3.2) is a consequence of (3.1)
by an argument similar to Proposition 4.3.1 in [KMN1].
In this paper we do not use this fact since we explicitly construct
the isomorphism for the cases we consider.
\vskip5mm

\vskip5mm
\noindent
{\bf Example 3.4.}
Let us consider the case of $\widehat{sl_3}$, $B_{\La_1}\ot B_{\La_2}$.
Here the energy function is given by the following.
$H=0$ on the elements
$$
	\normalbaselines\mmmth\offinterlineskip
	\vcenter{
   \hbox{$\Flect(0.7cm,0.7cm,\hbox{$\bullet$})$
         }\vskip-0.4pt
   \hbox{$\Flect(0.7cm,0.7cm,\hbox{$$})$
         }\vskip-0.4pt
   \hbox{$\Flect(0.7cm,0.7cm,\hbox{$$})$
         }
                }
\ot
	\vcenter{
   \hbox{$\Flect(0.7cm,0.7cm,\hbox{$\bullet$})$
         }\vskip-0.4pt
   \hbox{$\Flect(0.7cm,0.7cm,\hbox{$\bullet$})$
         }\vskip-0.4pt
   \hbox{$\Flect(0.7cm,0.7cm,\hbox{$$})$
         }
                },
\hskip1mm
	\vcenter{
   \hbox{$\Flect(0.7cm,0.7cm,\hbox{$\bullet$})$
         }\vskip-0.4pt
   \hbox{$\Flect(0.7cm,0.7cm,\hbox{$$})$
         }\vskip-0.4pt
   \hbox{$\Flect(0.7cm,0.7cm,\hbox{$$})$
         }
                }
\ot
	\vcenter{
   \hbox{$\Flect(0.7cm,0.7cm,\hbox{$\bullet$})$
         }\vskip-0.4pt
   \hbox{$\Flect(0.7cm,0.7cm,\hbox{$$})$
         }\vskip-0.4pt
   \hbox{$\Flect(0.7cm,0.7cm,\hbox{$\bullet$})$
         }
                },
\hskip1mm
	\vcenter{
   \hbox{$\Flect(0.7cm,0.7cm,\hbox{$$})$
         }\vskip-0.4pt
   \hbox{$\Flect(0.7cm,0.7cm,\hbox{$\bullet$})$
         }\vskip-0.4pt
   \hbox{$\Flect(0.7cm,0.7cm,\hbox{$$})$
         }
                }
\ot
	\vcenter{
   \hbox{$\Flect(0.7cm,0.7cm,\hbox{$\bullet$})$
         }\vskip-0.4pt
   \hbox{$\Flect(0.7cm,0.7cm,\hbox{$\bullet$})$
         }\vskip-0.4pt
   \hbox{$\Flect(0.7cm,0.7cm,\hbox{$$})$
         }
                },
\hskip1mm
	\vcenter{
   \hbox{$\Flect(0.7cm,0.7cm,\hbox{$$})$
         }\vskip-0.4pt
   \hbox{$\Flect(0.7cm,0.7cm,\hbox{$\bullet$})$
         }\vskip-0.4pt
   \hbox{$\Flect(0.7cm,0.7cm,\hbox{$$})$
         }
                }
\ot
	\vcenter{
   \hbox{$\Flect(0.7cm,0.7cm,\hbox{$\bullet$})$
         }\vskip-0.4pt
   \hbox{$\Flect(0.7cm,0.7cm,\hbox{$$})$
         }\vskip-0.4pt
   \hbox{$\Flect(0.7cm,0.7cm,\hbox{$\bullet$})$
         }
                },
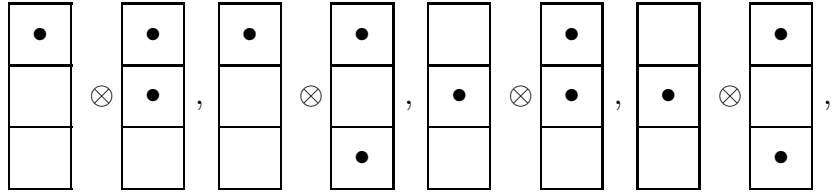
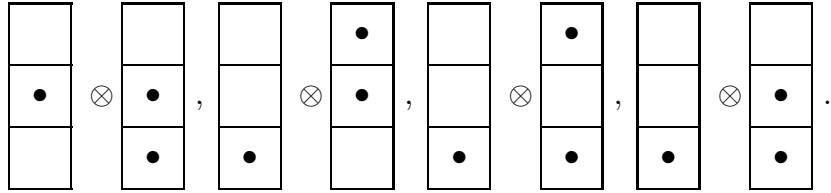
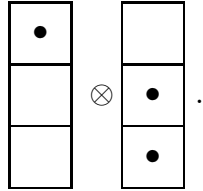
$$
and
$$
	\normalbaselines\mmmth\offinterlineskip
	\vcenter{
   \hbox{$\Flect(0.7cm,0.7cm,\hbox{$$})$
         }\vskip-0.4pt
   \hbox{$\Flect(0.7cm,0.7cm,\hbox{$\bullet$})$
         }\vskip-0.4pt
   \hbox{$\Flect(0.7cm,0.7cm,\hbox{$$})$
         }
                }
\ot
	\vcenter{
   \hbox{$\Flect(0.7cm,0.7cm,\hbox{$$})$
         }\vskip-0.4pt
   \hbox{$\Flect(0.7cm,0.7cm,\hbox{$\bullet$})$
         }\vskip-0.4pt
   \hbox{$\Flect(0.7cm,0.7cm,\hbox{$\bullet$})$
         }
                },
\hskip1mm
	\vcenter{
   \hbox{$\Flect(0.7cm,0.7cm,\hbox{$$})$
         }\vskip-0.4pt
   \hbox{$\Flect(0.7cm,0.7cm,\hbox{$$})$
         }\vskip-0.4pt
   \hbox{$\Flect(0.7cm,0.7cm,\hbox{$\bullet$})$
         }
                }
\ot
	\vcenter{
   \hbox{$\Flect(0.7cm,0.7cm,\hbox{$\bullet$})$
         }\vskip-0.4pt
   \hbox{$\Flect(0.7cm,0.7cm,\hbox{$\bullet$})$
         }\vskip-0.4pt
   \hbox{$\Flect(0.7cm,0.7cm,\hbox{$$})$
         }
                },
\hskip1mm
	\vcenter{
   \hbox{$\Flect(0.7cm,0.7cm,\hbox{$$})$
         }\vskip-0.4pt
   \hbox{$\Flect(0.7cm,0.7cm,\hbox{$$})$
         }\vskip-0.4pt
   \hbox{$\Flect(0.7cm,0.7cm,\hbox{$\bullet$})$
         }
                }
\ot
	\vcenter{
   \hbox{$\Flect(0.7cm,0.7cm,\hbox{$\bullet$})$
         }\vskip-0.4pt
   \hbox{$\Flect(0.7cm,0.7cm,\hbox{$$})$
         }\vskip-0.4pt
   \hbox{$\Flect(0.7cm,0.7cm,\hbox{$\bullet$})$
         }
                },
\hskip1mm
	\vcenter{
   \hbox{$\Flect(0.7cm,0.7cm,\hbox{$$})$
         }\vskip-0.4pt
   \hbox{$\Flect(0.7cm,0.7cm,\hbox{$$})$
         }\vskip-0.4pt
   \hbox{$\Flect(0.7cm,0.7cm,\hbox{$\bullet$})$
         }
                }
\ot
	\vcenter{
   \hbox{$\Flect(0.7cm,0.7cm,\hbox{$$})$
         }\vskip-0.4pt
   \hbox{$\Flect(0.7cm,0.7cm,\hbox{$\bullet$})$
         }\vskip-0.4pt
   \hbox{$\Flect(0.7cm,0.7cm,\hbox{$\bullet$})$
         }
                }.
$$
$H=1$ on
$$
	\normalbaselines\mmmth\offinterlineskip
	\vcenter{
   \hbox{$\Flect(0.7cm,0.7cm,\hbox{$\bullet$})$
         }\vskip-0.4pt
   \hbox{$\Flect(0.7cm,0.7cm,\hbox{$$})$
         }\vskip-0.4pt
   \hbox{$\Flect(0.7cm,0.7cm,\hbox{$$})$
         }
                }
\ot
	\vcenter{
   \hbox{$\Flect(0.7cm,0.7cm,\hbox{$$})$
         }\vskip-0.4pt
   \hbox{$\Flect(0.7cm,0.7cm,\hbox{$\bullet$})$
         }\vskip-0.4pt
   \hbox{$\Flect(0.7cm,0.7cm,\hbox{$\bullet$})$
         }
                }.
$$
The crystal graph of the tensor product is written in
{\bf Figure 3.1}. ($0$ arrows are omitted.)
The property of the energy function is obvious from this figure.

$$
\matrix{
{}&\bullet&\mapright1&\bullet&\mapright2&\bullet\cr
&&&&& \cr
\bullet&\bullet&\mapright1&\bullet&\mapright2&\bullet\cr
\downarrow2&\downarrow2&{}&{}&{}&\downarrow2\cr
\bullet&\bullet&\mapright1&\bullet&{}&\bullet\cr
\downarrow1&{}&{}&\downarrow1&{}&\downarrow1\cr
\bullet&\bullet&{}&\bullet&\mapright2&\bullet\cr
&&&&& \cr
{}&{}&{}&\hbox{{\bf figure 3.1}}&{}&{}\cr
}
$$
\vskip5mm

\noindent
{\bf 3.5. Explicit forms of the isomorphism $\iota$
and energy function $H$} \par
\noindent
First we demonstrate the connectivity for relevant cases.
\vskip5mm

\noindent
{\bf Proposition 3.6.}
$B_{\Lambda_k} \otimes B_{\Lambda_l}$ is connected. \par
\vskip3mm

\noindent
Proof.
Since $B_{\Lambda_k}$ and $B_{\Lambda_l}$ are both perfect crystals
of level one, $B_{\Lambda_k} \otimes B_{\Lambda_l}$ is connected
by Lemma 4.6.2 in [KMN1]. q.e.d.
\vskip5mm

\noindent
{\bf Proposition 3.7.}
$B_{k \Lambda_1} \otimes B_{l \Lambda_1}$ is connected. \par
\vskip3mm

\noindent
Proof. Since $B_{k \Lambda_1} \otimes B_{l \Lambda_1}$
is a union of highest weight crystals as a $\{1,\cdots,n-1\}$ crystal,
any element is connected to the element of the form
$$
	\normalbaselines\mmmth\offinterlineskip
	\vcenter{
   \hbox{$\Flect(0.7cm,0.7cm,\hbox{$1$})$\hskip-0.4pt
         $\Flect(0.7cm,1.4cm,\hbox{$\cdots$})$\hskip-0.4pt
         $\Flect(0.7cm,0.7cm,\hbox{$1$})$\hskip-0.4pt
         }
                }
\ot
	\vcenter{
   \hbox{$\Flect(0.7cm,0.7cm,\hbox{$i_1$})$\hskip-0.4pt
         $\Flect(0.7cm,1.4cm,\hbox{$\cdots$})$\hskip-0.4pt
         $\Flect(0.7cm,0.7cm,\hbox{$i_l$})$\hskip-0.4pt
         }
                }.
$$
By applying suitable $\ft_i$'s$(1\leq i\leq n-1)$
to this element, we arrive at the element of the form
$$
	\normalbaselines\mmmth\offinterlineskip
	\vcenter{
   \hbox{$\Flect(0.7cm,0.7cm,\hbox{$1$})$\hskip-0.4pt
         $\Flect(0.7cm,1.4cm,\hbox{$\cdots$})$\hskip-0.4pt
         $\Flect(0.7cm,0.7cm,\hbox{$1$})$\hskip-0.4pt
         }
                }
\ot
	\vcenter{
   \hbox{$\Flect(0.7cm,0.7cm,\hbox{$1$})$\hskip-0.4pt
         $\Flect(0.7cm,1.4cm,\hbox{$\cdots$})$\hskip-0.4pt
         $\Flect(0.7cm,0.7cm,\hbox{$1$})$\hskip-0.4pt
         $\Flect(0.7cm,0.7cm,\hbox{$n$})$\hskip-0.4pt
         $\Flect(0.7cm,1.4cm,\hbox{$\cdots$})$\hskip-0.4pt
         $\Flect(0.7cm,0.7cm,\hbox{$n$})$\hskip-0.4pt
         }
                }.
$$
Let $m$ be the number of $n$ in the second component.
Then if we apply $\ft_0^m$ to this element, we arrive at
$$
	\normalbaselines\mmmth\offinterlineskip
	\vcenter{
   \hbox{$\Flect(0.7cm,0.7cm,\hbox{$1$})$\hskip-0.4pt
         $\Flect(0.7cm,1.4cm,\hbox{$\cdots$})$\hskip-0.4pt
         $\Flect(0.7cm,0.7cm,\hbox{$1$})$\hskip-0.4pt
         }
                }
\ot
	\vcenter{
   \hbox{$\Flect(0.7cm,0.7cm,\hbox{$1$})$\hskip-0.4pt
         $\Flect(0.7cm,1.4cm,\hbox{$\cdots$})$\hskip-0.4pt
         $\Flect(0.7cm,0.7cm,\hbox{$1$})$\hskip-0.4pt
         }
                }.
$$
q.e.d.
\vskip5mm

By these propositions, if there exist isomorphisms
$$
\iota \ : \ B_{\Lambda_k} \otimes B_{\Lambda_l}
\rightarrow B_{\Lambda_l} \otimes B_{\Lambda_k},
$$
$$
\iota \ : \ B_{k \Lambda_1} \otimes B_{l \Lambda_1}
\rightarrow B_{l \Lambda_1} \otimes B_{k \Lambda_1},
$$
then the $\iota$'s are unique, and the corresponding energy function
is unique up to additive constants.
Now we shall give the explicit rules for the isomorphism
and the energy function.
We assume $k \geq l$.
The proofs for these rules are given in \S 3.15.
\vskip5mm

\noindent
{\bf Example 3.8.} \
The following are examples of the isomorphism
$$
\iota \ : \ B_{\Lambda_3} \otimes B_{\Lambda_2}
\rightarrow B_{\Lambda_2} \otimes B_{\Lambda_3}.
$$

$$
	\normalbaselines\mmmth\offinterlineskip
	\vcenter{
   \hbox{$\Flect(0.7cm,0.7cm,\hbox{$\bullet$})$
         }\vskip-0.4pt
   \hbox{$\Flect(0.7cm,0.7cm,\hbox{$$})$
         }\vskip-0.4pt
   \hbox{$\Flect(0.7cm,0.7cm,\hbox{$\bullet_a'$})$
         }\vskip-0.4pt
   \hbox{$\Flect(0.7cm,0.7cm,\hbox{$$})$
         }\vskip-0.4pt
   \hbox{$\Flect(0.7cm,0.7cm,\hbox{$\bullet_b'$})$
         }
                }
\otimes
	\vcenter{
   \hbox{$\Flect(0.7cm,0.7cm,\hbox{$$})$
         }\vskip-0.4pt
   \hbox{$\Flect(0.7cm,0.7cm,\hbox{$\bullet_a$})$
         }\vskip-0.4pt
   \hbox{$\Flect(0.7cm,0.7cm,\hbox{$\bullet_b$})$
         }\vskip-0.4pt
   \hbox{$\Flect(0.7cm,0.7cm,\hbox{$$})$
         }\vskip-0.4pt
   \hbox{$\Flect(0.7cm,0.7cm,\hbox{$$})$
         }
                },
\rightarrow
	\vcenter{
   \hbox{$\Flect(0.7cm,0.7cm,\hbox{$$})$
         }\vskip-0.4pt
   \hbox{$\Flect(0.7cm,0.7cm,\hbox{$$})$
         }\vskip-0.4pt
   \hbox{$\Flect(0.7cm,0.7cm,\hbox{$\bullet_a'$})$
         }\vskip-0.4pt
   \hbox{$\Flect(0.7cm,0.7cm,\hbox{$$})$
         }\vskip-0.4pt
   \hbox{$\Flect(0.7cm,0.7cm,\hbox{$\bullet_b'$})$
         }
                }
\otimes
	\vcenter{
   \hbox{$\Flect(0.7cm,0.7cm,\hbox{$\bullet$})$
         }\vskip-0.4pt
   \hbox{$\Flect(0.7cm,0.7cm,\hbox{$\bullet_a$})$
         }\vskip-0.4pt
   \hbox{$\Flect(0.7cm,0.7cm,\hbox{$\bullet_b$})$
         }\vskip-0.4pt
   \hbox{$\Flect(0.7cm,0.7cm,\hbox{$$})$
         }\vskip-0.4pt
   \hbox{$\Flect(0.7cm,0.7cm,\hbox{$$})$
         }
                },
\hskip1cm
	\vcenter{
   \hbox{$\Flect(0.7cm,0.7cm,\hbox{$\bullet_b'$})$
         }\vskip-0.4pt
   \hbox{$\Flect(0.7cm,0.7cm,\hbox{$\bullet$})$
         }\vskip-0.4pt
   \hbox{$\Flect(0.7cm,0.7cm,\hbox{$$})$
         }\vskip-0.4pt
   \hbox{$\Flect(0.7cm,0.7cm,\hbox{$\bullet_a'$})$
         }\vskip-0.4pt
   \hbox{$\Flect(0.7cm,0.7cm,\hbox{$$})$
         }
                }
\otimes
	\vcenter{
   \hbox{$\Flect(0.7cm,0.7cm,\hbox{$$})$
         }\vskip-0.4pt
   \hbox{$\Flect(0.7cm,0.7cm,\hbox{$$})$
         }\vskip-0.4pt
   \hbox{$\Flect(0.7cm,0.7cm,\hbox{$\bullet_a$})$
         }\vskip-0.4pt
   \hbox{$\Flect(0.7cm,0.7cm,\hbox{$$})$
         }\vskip-0.4pt
   \hbox{$\Flect(0.7cm,0.7cm,\hbox{$\bullet_b$})$
         }
                }
\rightarrow
	\vcenter{
   \hbox{$\Flect(0.7cm,0.7cm,\hbox{$\bullet_b'$})$
         }\vskip-0.4pt
   \hbox{$\Flect(0.7cm,0.7cm,\hbox{$$})$
         }\vskip-0.4pt
   \hbox{$\Flect(0.7cm,0.7cm,\hbox{$$})$
         }\vskip-0.4pt
   \hbox{$\Flect(0.7cm,0.7cm,\hbox{$\bullet_a'$})$
         }\vskip-0.4pt
   \hbox{$\Flect(0.7cm,0.7cm,\hbox{$$})$
         }
                }
\otimes
	\vcenter{
   \hbox{$\Flect(0.7cm,0.7cm,\hbox{$$})$
         }\vskip-0.4pt
   \hbox{$\Flect(0.7cm,0.7cm,\hbox{$\bullet$})$
         }\vskip-0.4pt
   \hbox{$\Flect(0.7cm,0.7cm,\hbox{$\bullet_a$})$
         }\vskip-0.4pt
   \hbox{$\Flect(0.7cm,0.7cm,\hbox{$$})$
         }\vskip-0.4pt
   \hbox{$\Flect(0.7cm,0.7cm,\hbox{$\bullet_b$})$
         }
                }.
$$
For the energy function, $H=0$ for the first example
and $H=-1$ for the second.
\vskip5mm

\noindent
{\bf Example 3.9.} \
The following are examples of the isomorphism
$$
\iota \ : \ B_{3 \Lambda_1} \otimes B_{2 \Lambda_1}
\rightarrow B_{2 \Lambda_1} \otimes B_{3 \Lambda_1}.
$$
$$
	\normalbaselines\mmmth\offinterlineskip
	\vcenter{
   \hbox{$\Flect(0.7cm,1.1cm,\hbox{$\bullet_a' \bullet$})$
         }\vskip-0.4pt
   \hbox{$\Flect(0.7cm,1.1cm,\hbox{$\bullet_b'$})$
         }\vskip-0.4pt
   \hbox{$\Flect(0.7cm,1.1cm,\hbox{$$})$
         }
                }
\otimes
	\vcenter{
   \hbox{$\Flect(0.7cm,1.1cm,\hbox{$$})$
         }\vskip-0.4pt
   \hbox{$\Flect(0.7cm,1.1cm,\hbox{$\bullet_a$})$
         }\vskip-0.4pt
   \hbox{$\Flect(0.7cm,1.1cm,\hbox{$\bullet_b$})$
         }
                }
\rightarrow
	\vcenter{
   \hbox{$\Flect(0.7cm,1.1cm,\hbox{$\bullet_a'$})$
         }\vskip-0.4pt
   \hbox{$\Flect(0.7cm,1.1cm,\hbox{$\bullet_b'$})$
         }\vskip-0.4pt
   \hbox{$\Flect(0.7cm,1.1cm,\hbox{$$})$
         }
                }
\otimes
	\vcenter{
   \hbox{$\Flect(0.7cm,1.1cm,\hbox{$\bullet$})$
         }\vskip-0.4pt
   \hbox{$\Flect(0.7cm,1.1cm,\hbox{$\bullet_a$})$
         }\vskip-0.4pt
   \hbox{$\Flect(0.7cm,1.1cm,\hbox{$\bullet_b$})$
         }
                },
\hskip1cm
	\vcenter{
   \hbox{$\Flect(0.7cm,1.1cm,\hbox{$$})$
         }\vskip-0.4pt
   \hbox{$\Flect(0.7cm,1.1cm,\hbox{$\bullet_b' \bullet$})$
         }\vskip-0.4pt
   \hbox{$\Flect(0.7cm,1.1cm,\hbox{$\bullet_a'$})$
         }
                }
\otimes
	\vcenter{
   \hbox{$\Flect(0.7cm,1.1cm,\hbox{$\bullet_a$})$
         }\vskip-0.4pt
   \hbox{$\Flect(0.7cm,1.1cm,\hbox{$\bullet_b$})$
         }\vskip-0.4pt
   \hbox{$\Flect(0.7cm,1.1cm,\hbox{$$})$
         }
                }
\rightarrow
	\vcenter{
   \hbox{$\Flect(0.7cm,1.1cm,\hbox{$$})$
         }\vskip-0.4pt
   \hbox{$\Flect(0.7cm,1.1cm,\hbox{$\bullet_b'$})$
         }\vskip-0.4pt
   \hbox{$\Flect(0.7cm,1.1cm,\hbox{$\bullet_a'$})$
         }
                }
\otimes
	\vcenter{
   \hbox{$\Flect(0.7cm,1.1cm,\hbox{$\bullet_a$})$
         }\vskip-0.4pt
   \hbox{$\Flect(0.7cm,1.1cm,\hbox{$\bullet \bullet_b$})$
         }\vskip-0.4pt
   \hbox{$\Flect(0.7cm,1.1cm,\hbox{$$})$
         }
                }.
$$
For the first example $H=0$, and for the second $H=2$.
\vskip5mm

\noindent
{\bf Rule 3.10.}
The case of $B_{\Lambda_k} \otimes B_{\Lambda_l}$ $(k\geq l)$.
\par
\noindent
The general procedure to obtain the isomorphism and energy function
is as follows.

{\bf 1.} \ Pick the dot $\bullet_a$
which has the highest position in the 2nd column and connect it
with the dot $\bullet_a'$ in the first column which
has the highest position
among all dots whose positions are not higher than $\bullet_a$.
If there is no such dot in the 1st column,
return to the top, assuming the periodic boundary condition.
We call such a pair a "winding".

{\bf 2.} \ Repeat the same procedure for the remaining unconnected dots
$l-1$-times.

{\bf 3.} \ The isomorphism $\iota$ is obtained by sliding
the remaining $(k-l)$ unpaired dots in the 1st column to
the 2nd.

{\bf 4.} \ The energy function is $(-1)$ times
the number of "winding" pairs.
\vskip5mm

\noindent
{\bf Rule 3.11.}
The case of $B_{k \Lambda_1} \otimes B_{l \Lambda_1}$ $(k\geq l)$.
\par

\noindent
In this case the rule to determine the isomorphism
$\iota$ and  energy function
$H$ is almost the same as that in the anti-symmetric case.
The only differences are in {\bf 1.} and {\bf 4}.

{\bf 1.} \
The partner $\bullet_a'$ is the
dot which has the lowest position among all dots
whose positions are higher than that of $\bullet_a$.
If there is no such dot return to the bottom.

{\bf 4.} \
The energy function is simply the number of "winding" pairs.
\vskip5mm

\noindent
{\bf Definition 3.12.} \
We refer to the line appearing in Rules 3.10 and 3.11
(connecting two dots) as the $H$-line.
\vskip2mm

\noindent
Any way of drawing the $H$-lines is allowed as long as
it represents Rules 3.10 and 3.11 correctly, however, we
sometimes use the following special way of drawing $H$-lines,
especially in the proof of Lemma 3.16 and in Section 5
for the sake of convenience.
This special way is as follows.
Starting from a dot in the 2nd column, first go horizontally left
and then go vertically down, until you find its partner in the first column.
Then, if necessary, round up outside of the diagram and again
go vertically down from the top.
\vskip5mm

\noindent
{\bf 3.13. The relation of energy functions}
\par
\noindent
We denote by $H_{\La_k\La_l}$ and $H_{k\La_1l\La_1}$ the energy
functions of
$B_{\La_k}\ot B_{\La_l}$ and $B_{k\La_1}\ot B_{l\La_1}$, respectively,
which are explicitly described above. They are characterized
by the normalization $H_{\La_k\La_l}=0$ on
$$
	\normalbaselines\mmmth\offinterlineskip
	\vcenter{
   \hbox{$\Flect(0.7cm,0.7cm,\hbox{$1$})$
         }\vskip-0.4pt
   \hbox{$\Flect(1.4cm,0.7cm,\hbox{$\vdots$})$
         }\vskip-0.4pt
   \hbox{$\Flect(0.7cm,0.7cm,\hbox{$k$})$
         }
                }
\ot
	\vcenter{
   \hbox{$\Flect(0.7cm,0.7cm,\hbox{$1$})$
         }\vskip-0.4pt
   \hbox{$\Flect(1.4cm,0.7cm,\hbox{$\vdots$})$
         }\vskip-0.4pt
   \hbox{$\Flect(0.7cm,0.7cm,\hbox{$l$})$
         }
                }
$$
and $H_{k\La_1l\La_1}=l$ on
$$
	\vcenter{
   \hbox{$\Flect(0.7cm,0.7cm,\hbox{$1$})$\hskip-0.4pt
         $\Flect(0.7cm,1.4cm,\hbox{$\cdots$})$\hskip-0.4pt
         $\Flect(0.7cm,0.7cm,\hbox{$1$})$
         }
                }
\ot
	\vcenter{
   \hbox{$\Flect(0.7cm,0.7cm,\hbox{$1$})$\hskip-0.4pt
         $\Flect(0.7cm,1.4cm,\hbox{$\cdots$})$\hskip-0.4pt
         $\Flect(0.7cm,0.7cm,\hbox{$1$})$
         }
                }.
$$
Consider the semi-standard tableau $T\in {\cal T}(\la,\mu)$, and
let $b=b_1\ot\cdots \ot b_r$ and
$\tilde{b}=\tilde{b}_1\ot\cdots \ot\tilde{b}_r$ be
the corresponding highest weight elements
in $B_{\La_{\mu_1}}\ot\cdots\ot B_{\La_{\mu_r}}$ and
$B_{\mu_1\La_{1}}\ot\cdots\ot B_{\mu_r\La_{r}}$, respectively.
We define $b_i^{(j)}$ and $\tilde{b}_i^{(j)}$ as in the
beginning of Section 4.
Then we have
\vskip5mm

\noindent
{\bf Proposition 3.14.}
The following relation holds:
$$
H_{\La_j\La_i}(b_j\ot b_{i}^{(j+1)})+
H_{j\La_1i\La_1}(\tilde{b}_j\ot \tilde{b}_{i}^{(j+1)})=0.
\eqno(3.3)
$$
\vskip3mm

\noindent
Proof. Comparing the evaluation rule of
$H_{\La_j\La_i}(b_j\ot b_{i}^{(j+1)})$
and $H_{j\La_1i\La_1}(\tilde{b}_j\ot \tilde{b}_{i}^{(j+1)})$
and the corresponding isomorphism rule on the semi-standard tableau
$T$, one finds that they are the same (see, for example, both rules
by Examples 2.10 and 2.11).
Hence beginning from $b$ and $\tilde{b}$, we have equation (3.3)
inductively. q.e.d.
\vskip5mm

\noindent
{\bf 3.15. Proof of Rules for $\iota$ and the energy function}
\par
\noindent
Here we prove that the rules described above give
the correct isomorphism $\iota$ and energy function.
In Rule 3.10 and 3.11, we made pairs of dots starting from the top.
However, this choice is not essential.
In fact we can again consider these rules by changing
the order in which the $H$ lines are drawn.
Then one can prove
\vskip5mm

\noindent
{\bf Lemma 3.16.} \
The isomorphism $\iota$ and the energy function
determined by the above rules are independent of
the order of drawing lines.
\vskip3mm

\noindent
Proof. \ (Anti-symmetric case.)\
Suppose there exists a dot (say $\bullet_{L1}$ ) in
one of the left boxes that is
unpaired in one ordering (say A) and paired in another (say B).
Let $\bullet_{R1}$ be the partner of $\bullet_{L1}$ in B, and
let $\bullet_{L2}$ be the partner of $\bullet_{R1}$ in A.
Then $\bullet_{L2}$ must be paired with some dot (say $\bullet_{R2}$ )
in B since the $R_1L_1$ line already passes through $\bullet_{L2}$
({\bf Figure 3.2}).
Let us prove that this process of determining $R_1L_2R_2\cdots$
does not stop. Suppose that we already have the $R_jL_j$ line in $B$.
Let us consider the following statement for $R_j$.
\vskip2mm
\noindent
$(Ass)_j$ \  For each box in the left column in $B$ whose position
is not higher than that of $R_j$ and not lower than that of $L_1$,
there is a line passing through it.(Not necessarily the end point.)
\par
\vskip2mm

\noindent
Here we assume periodic boundary conditions(pbc).
This implies, for example, that if $R_j$ is below $L_1$ in the
non-pbc sense, then the position specified by $(Ass)_j$ is like that
in {\bf Figure 3.3}.
Suppose that $(Ass)_j$ holds. Let us prove that there then exist $L_{j+1}$
and $R_{j+1}$ for which $(Ass)_{j+1}$ holds.
By definition $L_{j+1}$ exists. Note that the position
of $L_{j+1}$ is above that of $L_1$ and not higher than that of $R_j$,
since $L_1$ is not paired in $A$. Hence, in $B$, some line
should already pass through $L_{j+1}$ by $(Ass)_j$,
and $R_{j+1}$ is determined.
If the position of $L_{j+1}$ is above that of $L_1$ in the non-pbc sense,
then $(Ass)_{j+1}$ obviously holds by $(Ass)_j$ and the drawing rule
({\bf Figure 3.4}).
If the position of $L_{j+1}$ is lower than the position of
$L_1$ in the non-pbc sense,
again $(Ass)_{j+1}$ obviously holds ({\bf Figure 3.5-6}).
By definition $R_1R_2\cdots$ are all distinct. This is a contradiction,
since the number of boxes in the right column is finite.
Therefore the set of end points of $H$ lines is independent of the
order in which they are drawn. This proves that the rule determining $\iota$
is independent of the order of drawing lines.

Next let us prove that the value of the energy function $H$
determined by the above
rule is independent of the order of drawing the $H$ lines.

Let $l_0,\cdots,l_n$ be the horizontal lines in the diagram
numbered from the top down. To each $l_j$ we assign the non-negative
integer $h_J(l_j)$ as the number of $H$ lines crossing $l_j$
({\bf Figure 3.7}).
Here the subscript $J$ specifies an order.
The value of the energy function $H$ w.r.t. $J$ is given by $h_J(l_0)$.
By the rule of drawing $H$ lines, $h_J(l)$ is subject to the following
relation:
\vskip2mm
\itemitem{(Rel)} $h_J(l_{j+1})=h_J(l_j)+e^r_j-e^l_j$,
\par
\vskip2mm
\noindent
where $e^l_j$ is the number of $\bullet$ which are end points
of the $H$ lines and are between $l_j$ and $l_{j+1}$.
Similarly $e^r_j$ is the number of $\bullet$ between $l_j$ and $l_{j+1}$.
In particular $e^r_j,e^l_j=0,1$.

The set $(h_J(l_0),\cdots,h_J(l_{n-1}))$ is the same for any
$J$ up to an over all additive constant by (Rel), since,
as previously proven, the set of end points of $H$ lines
is independent of $J$.
We prove that $h_J(l_j)$ does not depend on $J$.
To this end it is sufficient to prove the existence of $l_j$ for
which $h_J(l_j)=0$.
In fact for two orders $J_1$ and $J_2$, if there is an $i$ such that
$
h_{J_1}(l_i)=h_{J_2}(l_i)+m,\ m>0,
$
then by the above comment,
$$
h_{J_1}(l_k)=h_{J_2}(l_k)+m\quad \hbox{ for any $k$}.
$$
Since $h_{J_2}(l_k)\geq0$ for any $k$, this contradicts the existence
of $l_j$ with $h_{J_1}(l_j)=0$.
So let us prove the existence of such $l_j$.

Let $\bullet$ (say $L$) be the end point (in the left column)
of a winding line and is lower than all other such points.
Let $\bullet$ (say $R$) be the point in the right column
from which a winding line starts and is higher than all other such points.
We claim that there exists a horizontal line $l$
between the underline of the box $L$ and the overline of the box $R$
for which $h_J(l)=0$.

Suppose that such an $l$ does not exist.
Then the situation is like that in {\bf Figure 3.8}.
We proceed with our proof according to this figure.
Note first that the order of $\alpha_2$ is smaller than
or equal to that of $\alpha_1$ by the definition of
$R$ and $L$.
Next note that
\vskip2mm
\item{$\bullet$}
the order of $\alpha_3$ is larger than that of $\alpha_1$.
\item{$\bullet$}
the order of $\alpha_4$ is larger than that of $\alpha_3$.
\item{$\bullet$}
the order of $\alpha_5$ is larger than that of $\alpha_4$.
\item{$\bullet$}
the order of $\alpha_2$ is larger than that of $\alpha_5$.
\par
\vskip2mm

\noindent
Hence the order of $\alpha_2$ is larger than $\alpha_1$, which
is a contradiction.
Thus in the anti-symmetric case, the value of the energy function
given by our rule does not depend on the order of drawing lines.
\vskip2mm

\noindent
(Symmetric case.)\
The proof of the symmetric case reduces to the anti-symmetric case
in the following way.
Let us make the following change in the diagram
(see {\bf Figure 3.9}):
\vskip2mm
\itemitem{1.}
Within each box, align the $\bullet$ vertically.
Then, shift the left column down by one box.
Create a rectangular diagram from the present one by adding two boxes.
\itemitem{2.} Draw horizontal lines such that each $\bullet$
is in one box.
\itemitem{3.} Reverse up and down in the diagram.
\par
\vskip2mm

\noindent
Now for the final diagram, the drawing rule for $H$ lines
is the same as that in the anti-symmetric case.
Hence the independence of the value of the $H$ function
and the isomorphism $\iota$ follows from the anti-symmetric case.
q.e.d.
\vskip5mm

\noindent
{\bf Proposition 3.17.}
The map $\iota$ described in the above rules gives the isomorphism
of crystals, i.e., it commutes with the actions of ${\tilde {f_i}}$
and  ${\tilde {e_i}}$.
\vskip3mm

\noindent
Proof. (Anti-symmetric case.)
In the following proof (also for the symmetric case)
we prove the commutativity with $\ft_i$.
The case for $\et_i$ is similarly proved.
There are four cases of the
configurations of the $i$-th boxes,
$$
	\normalbaselines\mmmth\offinterlineskip
	\vcenter{
   \hbox{$\Flect(0.7cm,0.7cm,\hbox{$$})$
         }
                }
\otimes
	\vcenter{
   \hbox{$\Flect(0.7cm,0.7cm,\hbox{$$})$
         }
                },
\quad
	\vcenter{
   \hbox{$\Flect(0.7cm,0.7cm,\hbox{$$})$
         }
                }
\otimes
	\vcenter{
   \hbox{$\Flect(0.7cm,0.7cm,\hbox{$\bullet$})$
         }
                },
\quad
	\vcenter{
   \hbox{$\Flect(0.7cm,0.7cm,\hbox{$\bullet$})$
         }
                }
\otimes
	\vcenter{
   \hbox{$\Flect(0.7cm,0.7cm,\hbox{$$})$
         }
                },
\quad
	\vcenter{
   \hbox{$\Flect(0.7cm,0.7cm,\hbox{$\bullet$})$
         }
                }
\otimes
	\vcenter{
   \hbox{$\Flect(0.7cm,0.7cm,\hbox{$\bullet$})$
         }
                }.
$$
The first case is obvious (${\tilde {f_i}}$=0).
For the 2nd case, it is easy to see that
the $i$-th and $(i+1)$-th configuration of $b\ot b^\prime$ is invariant
under $\iota$ regardless of the $(i+1)$-th configuration in $b \otimes b'$.
Hence $\iota$ commutes with ${\tilde {f_i}}$.
For the 3rd case, possible nontrivial patterns (i.e.,
those in which the $i$-th and $(i+1)$-th configurations change under $\iota$)
are
$$
	\normalbaselines\mmmth\offinterlineskip
	\vcenter{
   \hbox{$\Flect(0.7cm,0.7cm,\hbox{$\bullet$})$
         }\vskip-0.4pt
   \hbox{$\Flect(0.7cm,0.7cm,\hbox{$$})$
         }
                }
\otimes
	\vcenter{
   \hbox{$\Flect(0.7cm,0.7cm,\hbox{$$})$
         }\vskip-0.4pt
   \hbox{$\Flect(0.7cm,0.7cm,\hbox{$$})$
         }
                }
\rightarrow
	\vcenter{
   \hbox{$\Flect(0.7cm,0.7cm,\hbox{$$})$
         }\vskip-0.4pt
   \hbox{$\Flect(0.7cm,0.7cm,\hbox{$$})$
         }
                }
\otimes
	\vcenter{
   \hbox{$\Flect(0.7cm,0.7cm,\hbox{$\bullet$})$
         }\vskip-0.4pt
   \hbox{$\Flect(0.7cm,0.7cm,\hbox{$$})$
         }
                },
\hskip1cm
	\vcenter{
   \hbox{$\Flect(0.7cm,0.7cm,\hbox{$\bullet$})$
         }\vskip-0.4pt
   \hbox{$\Flect(0.7cm,0.7cm,\hbox{$\bullet$})$
         }
                }
\otimes
	\vcenter{
   \hbox{$\Flect(0.7cm,0.7cm,\hbox{$$})$
         }\vskip-0.4pt
   \hbox{$\Flect(0.7cm,0.7cm,\hbox{$\bullet$})$
         }
                }
\rightarrow
	\vcenter{
   \hbox{$\Flect(0.7cm,0.7cm,\hbox{$$})$
         }\vskip-0.4pt
   \hbox{$\Flect(0.7cm,0.7cm,\hbox{$\bullet$})$
         }
                }
\otimes
	\vcenter{
   \hbox{$\Flect(0.7cm,0.7cm,\hbox{$\bullet$})$
         }\vskip-0.4pt
   \hbox{$\Flect(0.7cm,0.7cm,\hbox{$\bullet$})$
         }
                },
$$
and
$$
	\normalbaselines\mmmth\offinterlineskip
	\vcenter{
   \hbox{$\Flect(0.7cm,0.7cm,\hbox{$\bullet$})$
         }\vskip-0.4pt
   \hbox{$\Flect(0.7cm,0.7cm,\hbox{$\bullet$})$
         }
                }
\otimes
	\vcenter{
   \hbox{$\Flect(0.7cm,0.7cm,\hbox{$$})$
         }\vskip-0.4pt
   \hbox{$\Flect(0.7cm,0.7cm,\hbox{$$})$
         }
                }
\rightarrow
	\vcenter{
   \hbox{$\Flect(0.7cm,0.7cm,\hbox{$\bullet$})$
         }\vskip-0.4pt
   \hbox{$\Flect(0.7cm,0.7cm,\hbox{$$})$
         }
                }
\otimes
	\vcenter{
   \hbox{$\Flect(0.7cm,0.7cm,\hbox{$$})$
         }\vskip-0.4pt
   \hbox{$\Flect(0.7cm,0.7cm,\hbox{$\bullet$})$
         }
                }
\quad {\rm or} \quad
	\vcenter{
   \hbox{$\Flect(0.7cm,0.7cm,\hbox{$$})$
         }\vskip-0.4pt
   \hbox{$\Flect(0.7cm,0.7cm,\hbox{$$})$
         }
                }
\otimes
	\vcenter{
   \hbox{$\Flect(0.7cm,0.7cm,\hbox{$\bullet$})$
         }\vskip-0.4pt
   \hbox{$\Flect(0.7cm,0.7cm,\hbox{$\bullet$})$
         }
                },
\hskip1cm
	\vcenter{
   \hbox{$\Flect(0.7cm,0.7cm,\hbox{$\bullet$})$
         }\vskip-0.4pt
   \hbox{$\Flect(0.7cm,0.7cm,\hbox{$$})$
         }
                }
\otimes
	\vcenter{
   \hbox{$\Flect(0.7cm,0.7cm,\hbox{$$})$
         }\vskip-0.4pt
   \hbox{$\Flect(0.7cm,0.7cm,\hbox{$\bullet$})$
         }
                }
\rightarrow
	\vcenter{
   \hbox{$\Flect(0.7cm,0.7cm,\hbox{$$})$
         }\vskip-0.4pt
   \hbox{$\Flect(0.7cm,0.7cm,\hbox{$$})$
         }
                }
\otimes
	\vcenter{
   \hbox{$\Flect(0.7cm,0.7cm,\hbox{$\bullet$})$
         }\vskip-0.4pt
   \hbox{$\Flect(0.7cm,0.7cm,\hbox{$\bullet$})$
         }
                }.
$$
One can easily check the commutativity for all cases.
For the last case, there is only one nontrivial pattern :
$$
	\normalbaselines\mmmth\offinterlineskip
	\vcenter{
   \hbox{$\Flect(0.7cm,0.7cm,\hbox{$\bullet$})$
         }\vskip-0.4pt
   \hbox{$\Flect(0.7cm,0.7cm,\hbox{$\bullet$})$
         }
                }
\otimes
	\vcenter{
   \hbox{$\Flect(0.7cm,0.7cm,\hbox{$\bullet$})$
         }\vskip-0.4pt
   \hbox{$\Flect(0.7cm,0.7cm,\hbox{$$})$
         }
                }
\rightarrow
	\vcenter{
   \hbox{$\Flect(0.7cm,0.7cm,\hbox{$\bullet$})$
         }\vskip-0.4pt
   \hbox{$\Flect(0.7cm,0.7cm,\hbox{$$})$
         }
                }
\otimes
	\vcenter{
   \hbox{$\Flect(0.7cm,0.7cm,\hbox{$\bullet$})$
         }\vskip-0.4pt
   \hbox{$\Flect(0.7cm,0.7cm,\hbox{$\bullet$})$
         }
                }
$$
which also commutes with ${\tilde {f_i}}$.
\vskip3mm

\noindent
(Symmetric case.)
For $b \otimes b'$, The action of ${\tilde {f_i}}$ is determined
by the number of dots in $i$-th boxes that are not connected
to $(i+1)$-th dots. If we note that
the rule for $\iota$ preserves these numbers,
the commutativity with ${\tilde {f_i}}$ is easily proved.
q.e.d.
\vskip5mm

\noindent
{\bf Proposition 3.18.}
The value of the energy function $H$ given by Rules 3.10 and 3.11
has the correct properties of an energy function.
\vskip3mm
\noindent
Proof. (Anti-symmetric case.)
The invariance of the value of $H$ under the non-trivial (i.e. non-zero)
action of $\et_i$ and $\ft_i$($i\neq0$) can easily be checked.
Let us demonstrate the property of $H$ under the action of $\ft_0$.
There are three patterns (say, $RR$, $LR$ and $LL$)
for the nontrivial action of ${\tilde {f_0}}$
on $b \otimes b'$ and $\iota (b \otimes b')$,
where $LR$ means ${\tilde {f_0}}$ acts on the left component
of $b \otimes b'$ and the right component of $\iota (b \otimes b')$.

The property to be proved is
$$
H({\tilde {f_0}} (b \otimes b'))-H(b \otimes b')
=\left\{
\eqalign{
1 & \quad ({\rm for} \quad RR), \cr
0 & \quad ({\rm for} \quad LR), \cr
-1 & \quad ({\rm for} \quad LL). \cr
}
\right.
$$
This can be checked case by case as in {\bf Figure 3.10}.
In the figure, the boxes are the $n$($\equiv 0$)-th and $1$st ones.
\vskip3mm

\noindent
(Symmetric case.)
In this case the invariance under the action of
$\et_i$ and $\ft_i$($i\neq0$) can also be easily checked.
Here, there are also three patterns of the ${\tilde {f_0}}$-action.
The property to be proved is the same as that above.
This can also be checked case by case, as in {\bf Figure 3.11}.
q.e.d.
\vskip5mm

Sometimes it is useful to represent the isomorphism
$\iota \ : \ b_1 \otimes b_2 \rightarrow b_2' \otimes b_1'$
and the value of the energy function $H(b_1,b_2)=h$ by the following
diagram
$$
\matrix{
& \hfil b_1' \hfil & \cr
& \hfil \vert \hfil \cr
\hfil b_2' \hfil & \hfil {-h-} \hfil & \hfil b_2 \hfil \cr
& \hfil \vert \hfil \cr
& \hfil b_1 \hfil & \cr}
$$

For a sequence of crystals $B_1,\cdots, B_i$ the isomorphism
$$
\iota \ : \ (B_1 \otimes B_2 \otimes \cdots \otimes B_{i-1}) \otimes B_i
\rightarrow B_i \otimes (B_1 \otimes B_2 \otimes \cdots \otimes B_{i-1})
$$
is defined by the composition of adjacent isomorphisms as
$$
\matrix{
&\hfil b_1' \hfil&&\hfil b_2' \hfil&&\hfil b_{i-1}' \hfil \cr
&\hfil \vert \hfil&&\hfil \vert \hfil&&\hfil \vert \hfil \cr
\hfil b_i^{(1)} \hfil&\hfil {-h_1-} \hfil&
\hfil b_i^{(2)} \hfil&\hfil {-h_2-} \hfil&\cdots&
\hfil {-h_{i-1}-} \hfil&\hfil b_i=b_i^{(i)} \hfil \cr
&\hfil \vert \hfil&&\hfil \vert \hfil&&\hfil \vert \hfil \cr
&\hfil b_1 \hfil&&\hfil b_2 \hfil&&\hfil b_{i-1} \hfil \cr
}.
$$

This $L$-operator like construction plays an essential role in our
representation of the index in terms of the energy function. That is,
if the sequence $b_1 \otimes b_2 \otimes \cdots \in
B_1 \otimes B_2 \otimes \cdots $ is the highest weight element
corresponding to a given semi-standard tableau $T$,
the index of $i$-th component ind($i$) is given by
$$
{\rm ind}(i)=\pm (h_1+h_2+\cdots+h_{i-1}),
$$
where $+$ [resp. $-$] is for the symmetric [resp. anti-symmetric] case.

\vskip5mm
\noindent
{\bf Example 3.19.}

The following diagram indicates the way of calculating the ind($4$)
for the highest weight element in Example 2.11.
$$
\matrix{
&\hfil (134) \hfil&&\hfil (123) \hfil&&\hfill (45) \hfil \cr
&\hfil \vert \hfil&&\hfil \vert \hfil&&\hfil \vert \hfil \cr
\hfil (12) \hfil&\hfil {-1-} \hfil&\hfill (14) \hfil&
\hfil {-0-} \hfil&\hfil (13) \hfil&\hfil {-2-} \hfil&
\hfil (45) \hfil \cr
&\hfil \vert \hfil&&\hfil \vert \hfil&&\hfil \vert \hfil \cr
&\hfil (123) \hfil&&\hfil (124) \hfil&&\hfill (13) \hfil \cr
}
$$
where $(13)$ represents
$$
	\normalbaselines\mmmth\offinterlineskip
	\vcenter{
   \hbox{$\Flect(0.7cm,0.7cm,\hbox{$\bullet$})$
         }\vskip-0.4pt
   \hbox{$\Flect(0.7cm,0.7cm,\hbox{$$})$
         }\vskip-0.4pt
   \hbox{$\Flect(0.7cm,0.7cm,\hbox{$\bullet$})$
         }\vskip-0.4pt
   \hbox{$\Flect(0.7cm,0.7cm,\hbox{$$})$
         }\vskip-0.4pt
   \hbox{$\Flect(0.7cm,0.7cm,\hbox{$$})$
         }
                }
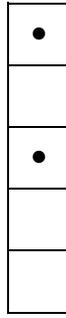
$$

\vfill\eject
%%%%%%%%%%%%%%%%%%%%%%%%%%%%%%%%%%

%Nakayashiki personal macros
\def\ot{\otimes}
\def\ra{\longrightarrow}
\def\et{\tilde{e}}
\def\ft{\tilde{f}}
\def\La{\Lambda}
\def\la{\lambda}
\def\mapright#1{\smash{\mathop{\longrightarrow}\limits^{#1}}}
\def\ind{\hbox{ind}}

\noindent
{\fontB 4. Index in Terms of the Energy Function}
\vskip2mm
\noindent
Let $T$ be a semistandard
tableau in ${\cal T}(\lambda,\mu)$ of shape $\lambda$ and weight
$\mu=(\mu_1,\cdots,\mu_l)$.
Let $b_1\ot\cdots\ot b_l$ in $B_{\La_{\mu_1}}\ot\cdots \ot B_{\La_{\mu_l}}$
and $\tilde{b}_1\ot\cdots\ot \tilde{b}_l$ in
$B_{\mu_1\La_1}\ot\cdots \ot B_{\mu_l\La_1}$
be the corresponding highest weight elements.
We define $b_{i}^{(j)}$ and $\tilde{b}_{i}^{(j)}$ by

$$
\matrix{
B_{\La_{\mu_1}}\ot\cdots \ot B_{\La_{\mu_i}}
\!\!&\simeq&\!\!
B_{\La_{\mu_1}}\ot\cdots \ot B_{\La_{\mu_i}}\ot B_{\La_{\mu_{i-1}}}
\!\!&\simeq&\!\!\cdots\!\!&\simeq&\!\!
B_{\La_{\mu_i}}\ot\cdots\ot B_{\La_{\mu_{i-1}}}
\cr
b_1\ot\cdots\ot b_i
\!\!&\mapsto&\!\!
b_1\ot\cdots\ot b_i^{(i-1)}\ot b_{i-1}^\prime
\!\!&\mapsto&\!\!\cdots\!\!&\mapsto&\!\!
b_i^{(1)}\ot b_1^\prime\ot\cdots\ot b_{i-1}^\prime
\cr
B_{\mu_1\La_{1}}\ot\cdots \ot B_{\mu_i\La_{1}}
\!\!&\simeq&\!\!
B_{\mu_1\La_{1}}\ot\cdots \ot B_{\mu_i\La_{1}}\ot B_{\mu_{i-1}\La_{1}}
\!\!&\simeq&\!\!\cdots\!\!&\simeq&\!\!
B_{\mu_i\La_{1}}\ot\cdots\ot B_{\mu_{i-1}\La_{1}}
\cr
\tilde{b}_1\ot\cdots\ot \tilde{b}_i
\!\!&\mapsto&\!\!
\tilde{b}_1\ot\cdots\ot \tilde{b}_i^{(i-1)}\ot \tilde{b}_{i-1}^\prime
\!\!&\mapsto&\!\!\cdots\!\!&\mapsto&\!\!
\tilde{b}_i^{(1)}\ot\tilde{b}_1^\prime\ot\cdots\ot \tilde{b}_{i-1}^\prime.
}
$$
We set $b_i^{(i)}=b_i$.
Then our main theorem is
\vskip5mm

\noindent
{\bf Theorem 4.1.}
For $1\leq i\leq l$ the following relation holds:
$$
\eqalign{
\ind(i)&=-\sum_{j=1}^{i-1}
H_{\La_{\mu_j}\La_{\mu_i}}(b_j,b_{i}^{(j+1)}),
\cr
&=\sum_{j=1}^{i-1}
H_{\mu_j\La_{1}\mu_i\La_{1}}(\tilde{b}_j,\tilde{b}_{i}^{(j+1)}).
}
$$
\vskip5mm

\noindent
The proof of this theorem is given in the next section.
As a consequence of this theorem and Theorem 2.16 of
Lascoux-Schutzenberger (and [KR,B] for the skew case)
we have
\vskip5mm

\noindent
{\bf Corollary 4.2.} There exist the following two expressions
of $K_{\lambda,\mu}(q)$:
$$
\eqalign{
K_{\lambda,\mu}(q)&=\sum_{b_1\ot\cdots\ot b_l\in {\cal T}(\lambda,\mu)}
q^{-\sum_{i=1}^l\sum_{j=1}^{i-1}
H_{\La_{\mu_j}\La_{\mu_i}}(b_j,b_{i}^{(j+1)})},
\cr
&=\sum_{\tilde{b}_1\ot\cdots\ot \tilde{b}_l\in {\cal T}(\lambda,\mu)}
q^{\sum_{i=1}^l\sum_{j=1}^{i-1}
H_{\mu_j\La_{1}\mu_i\La_{1}}(\tilde{b}_j,\tilde{b}_{i}^{(j+1)})}.
}
$$
More generally,
for the skew Kostka polynomial $K_{\lambda/\nu,\mu}(q)$ $(\nu\subset\la)$,
we have
$$
\eqalign{
K_{\lambda/\nu,\mu}(q)&
=\sum_{b_1\ot\cdots\ot b_l\in {\cal T}(\lambda/\nu,\mu)}
q^{-\sum_{i=1}^l\sum_{j=1}^{i-1}
H_{\La_{\mu_j}\La_{\mu_i}}(b_j,b_{i}^{(j+1)})},
\cr
&=\sum_{\tilde{b}_1\ot\cdots\ot \tilde{b}_l\in {\cal T}(\lambda/\nu,\mu)}
q^{\sum_{i=1}^l\sum_{j=1}^{i-1}
H_{\mu_j\La_{1}\mu_i\La_{1}}(\tilde{b}_j,\tilde{b}_{i}^{(j+1)})}.
}
$$

\vskip5mm

\noindent
As an immediate corollary of Corollary 4.2
we have
\vskip5mm

\noindent
{\bf Corollary 4.3.} For a dominant integral weight $\lambda$ of $sl_n$
let us set
$$
V(k\La_i:\lambda)=\{v\in V(k\La_i)\ \vert \ e_jv=0( 1\leq j\leq n-1),
\ \hbox{wt}v=\lambda\}.
$$
We assume $\vert\lambda\vert\equiv k i \hbox{ mod }n$.
Let $\la^{(N)}=
\big(\la_1+kN-{1\over n}(\vert\lambda\vert-k i),\cdots,
\la_n+kN-{1\over n}(\vert\lambda\vert-k i)\big)$
and
$\mu^{(N)}=\big(k^{nN+i}\big)$,
where $\la=(\la_1,\cdots,\la_n)$ $(\la_1\geq\cdots\geq\la_n\geq0)$.
Then we have
$$
\hbox{tr}_{V(k\La_i:\lambda)}(q^{-d})
=\lim_{N\ra\infty}q^{-A_N}K_{\la^{(N)}\mu^{(N)}}(q),
$$
where
$$
A_N={1\over2}knN(N-1)+kNi.
$$
\vskip5mm

The case $i=0$ (Theorem 1.1) was conjectured by A.N. Kirillov (Conjecture 4 in
[Kir]).
For the level one modules we have another expression
of the branching function as a limit of Kostka polynomials.
\vskip5mm

\noindent
{\bf Corollary 4.4.} Let $\la=(\la_1,\cdots,\la_l)$ be a partition
such that $\la_1\leq n$ and $\vert\lambda\vert\equiv0\hbox{ mod }n$.
Set $\lambda^{(N)}=\big(n,\cdots,n,\la_1,\cdots,\la_l\big)$,
$\mu^{(N)}=\big(k^{nN}\big)$, where the $n$ in front of $\la_1$
appears $kN-\vert\lambda\vert$ times.
Then for any $1\leq k\leq n-1$ we have
$$
\hbox{tr}_{V(\La_0:\lambda^\prime)}(q^{-d})
=\lim_{N\ra\infty}q^{-A^k_N}K_{\la^{(N)}\mu^{(N)}}(q^{-1}),
$$
where $\la^\prime$ is the transpose of $\la$, and
$$
A^k_N=\sum_{j=1}^{nN-1}jH_{\La_k\La_k}(p^{gr(k)}_{j+1},p^{gr(k)}_j).
$$
The path $p^{gr(k)}=(p^{gr(k)}_j)$ is the ground state path
which we explicitly describe below.
\vskip5mm

\noindent
In the path realization of $B(\La_0)$ by $B_{\La_k}$,
the ground state path $p^{gr(k)}=(p^{gr(k)}_j)$,
$p^{gr(k)}_j\in B_{\La_k}$ is
given explicitly by
$$
	\normalbaselines\mmmth\offinterlineskip
p^{gr(k)}_j=
	\vcenter{
   \hbox{$\Flect(0.7cm,2.0cm,\hbox{$(n-j)k+1$})$
         }\vskip-0.4pt
   \hbox{$\Flect(1.4cm,2.0cm,\hbox{$\vdots$})$
         }\vskip-0.4pt
   \hbox{$\Flect(0.7cm,2.0cm,\hbox{$(n-j+1)k$})$
         }
                }.
$$
To understand this correctly we set the following rules.
\vskip3mm
\item{(1)} The numbers in boxes are considered modulo $n$.
\item{(2)} We impose periodic boundary conditions.
\par
\vskip3mm
\noindent
By these rules one can consider $p^{gr(k)}_j$ as a
uniquely determined semi-standard tableau,
an element of $B_{\La_k}$.
\vskip5mm

\noindent
{\bf 4.5. Proof of Corollary 4.3}
\par
\noindent
Let us consider the case $\mu_j=k$ in the second
equation for $K_{\lambda,\mu}(q)$ in Corollary 4.2.
In this case $\tilde{b}_i^{(j)}=\tilde{b}_j$ for any $i$ and $j$.
Hence we have
$$
K_{\lambda,\mu}(q)
=\sum_{\tilde{b}_l\ot\cdots\ot \tilde{b}_1\in {\cal T}(\lambda,\mu)}
q^{\sum_{j=1}^ljH_{k\La_{1}k\La_{1}}(\tilde{b}_{j+1},\tilde{b}_j)}.
\eqno(4.1).
$$
Recall that $\tilde{b}_l\ot\cdots\ot \tilde{b}_1\in {\cal T}(\lambda,\mu)$
if and only if
$\tilde{b}_l\ot\cdots\ot \tilde{b}_1$ is a highest weight element
with weight $\la$ as an $sl_n$ crystal. Therefore
$$
\eqalign{
K_{\lambda^{(N)},\mu^{(N)}}(q)
&=\sum_{b_{nN+i}\ot\cdots\ot b_1\in {\cal H}(N)}
q^{\sum_{j=1}^{nN+i-1}jH_{k\La_{1}k\La_{1}}(b_{j+1},b_j)},
\cr
{\cal H}(N)&=
\{b\in B_{k\La_1}^{\ot(nN+i)}\vert \et_j(b)=0\,(1\leq j\leq n-1),
wtb=\la\}.
}
\eqno(4.2)
$$
We note that $B(k\La_i)\simeq B(k\La_0)\ot B_{k\La_1}^{\ot(nN+i)}$
by the path realization. Using this isomorphism and the definition
of a path, if we set
$$
{\cal H}(\infty)=
\{b\in B(k\La_i)\vert \et_j(b)=0\,(1\leq j\leq n-1),
wtb=\la\},
$$
we have
$$
\eqalign{
{\cal H}(\infty)&=\cup_{N=1}^\infty{\cal H}(N),
\cr
{\cal H}(N)&\subset{\cal H}(N+1).
}
$$
The embedding of ${\cal H}(N)$ into ${\cal H}(\infty)$
is given by sending $b\in {\cal H}(N)$ to
$b_{k\La_0}\ot b\in {\cal H}(\infty)$.
Now by the results of [JMMO,KMN1],
$$
\hbox{tr}_{V(k\La_i:\lambda)}(q^{-d})
=\sum_{b\in {\cal H}(\infty)}q^{\omega(b)}
=\lim_{N\ra\infty}\sum_{b\in {\cal H}(N)}q^{\omega_N(b)},
$$
where
$$
\eqalign{
\omega(b)&=
\sum_{j=1}^\infty
j\big(
H_{k\La_1k\La_1}(b_{j+1},b_j)-H_{k\La_1k\La_1}(p^{gr}_{j+1},p^{gr}_j)
\big),
\cr
\omega_N(b)&=
\sum_{j=1}^{nN+i}
j\big(
H_{k\La_1k\La_1}(b_{j+1},b_j)-H_{k\La_1k\La_1}(p^{gr}_{j+1},p^{gr}_j)
\big),
}\eqno(4.3)
$$
$b=(b_j)_{j=1}^\infty$ is the path realization and
$p^{gr}=(p^{gr}_j)_{j=1}^\infty$ is the ground state path given by
$$
\eqalign{
p^{gr}_j&=(i+1-j)^k,
\cr
j^k&=	\vcenter{
   \hbox{$\Flect(0.7cm,0.7cm,\hbox{$j$})$\hskip-0.4pt
         $\Flect(0.7cm,1.4cm,\hbox{$\cdots$})$\hskip-0.4pt
         $\Flect(0.7cm,0.7cm,\hbox{$j$})$\hskip-0.4pt
         }
                } \ .
}
$$
We consider $i+1-j$ modulo $n$. Note that if
$b_{nN+i}\ot\cdots\ot b_1\in{\cal H}(N)$, then $b_{nN+i}=1^k=p^{gr}_{nN+i}$
by the highest weight condition.
Note also that if $b\in{\cal H}(N)$ corresponds to
$(b_j)_{j=1}^\infty\in{\cal H}(\infty)$,
then $b_j=p^{gr}_j$ for all $j\geq nN+i+1$ by definition.
Thus if we set
$$
A_N=\sum_{j=1}^{nN+i-1}jH_{k\La_1k\La_1}(p^{gr}_{j+1},p^{gr}_j),
$$
we have from (4.2) and (4.3)
$$
q^{-A_N}K_{\lambda^{(N)},\mu^{(N)}}(q)
=\sum_{b\in{\cal H}(N)}q^{\omega_N(b)}.
$$
We can evaluate $A_N$ explicitly using
$$
H_{k\La_1k\La_1}((j-1)^k,j^k)=k\delta_{j1}
$$
as
$$
A_N={1\over2}knN(N-1)+kNi.
$$
q.e.d.
\vskip5mm

\noindent
{\bf 4.6. Proof of Corollary 4.4}
\par
\noindent
The proof of Corollary 4.4 is similar to that of Corollary 4.3.
In this case we use the first expression of $K_{\lambda,\mu}(q)$
in Corollary 4.2.
Then we have
$$
\eqalign{
K_{\lambda^{(N)},\mu^{(N)}}(q)
&=\sum_{b_{nN}\ot\cdots\ot b_1\in {\cal H}(N)}
q^{\sum_{j=1}^{nN-1}jH_{\La_{k}\La_{k}}(b_{j+1},b_j)},
\cr
{\cal H}(N)&=
\{b\in B_{k\La_k}^{\ot nN}\vert \et_j(b)=0\,(1\leq j\leq n-1),
wtb=\la^\prime\}.
}
$$
This time we use the path realization of $B(\La_0)$ by $B_{\La_k}$.
Then
$$
B(\La_0)=\{b=(b_j)_{j=1}^\infty\vert\, b_j\in B_{\La_k},\,
b_j=p^{gr(k)}_j\, (j>>0)\}.
$$
Define
$$
{\cal H}(\infty)=
\{b\in B(k\La_i)\vert \et_j(b)=0\,(1\leq j\leq n-1),
wtb=\la^\prime\}.
$$
Again we have the injection ${\cal H}(N)\ra {\cal H}(\infty)$
sending $b$ to $b_{\La_0}\ot b$.
Using the result [KMN1]
$$
\hbox{tr}_{V(\La_0:\lambda)}(q^{-d})
=\sum_{b\in {\cal H}(\infty)}q^{\omega(b)},
$$
we have Corollary 4.4. Here $\omega(b)$ is similarly defined
by replacing $H_{k\La_1k\La_1}(\cdot,\cdot)$
by $H_{\La_k\La_k}(\cdot,\cdot)$ in (4.3).
q.e.d.
%\bye

\vfill\eject

%Nakayashiki personal macros
\def\ot{\otimes}
\def\ra{\longrightarrow}
\def\et{\tilde{e}}
\def\ft{\tilde{f}}
\def\La{\Lambda}
\def\la{\lambda}
\def\mapright#1{\smash{\mathop{\longrightarrow}\limits^{#1}}}
\def\ind{\hbox{ind}}

\noindent
{\fontB 5. Proof of Theorem 4.1}
\vskip2mm
\noindent
In this section we give a proof of Theorem 4.1.
We first reduce the statement of the theorem to a
more tractable one.
\vskip5mm

\noindent
{\bf 5.1. Local index}

Let us introduce the local version of the index.
\vskip5mm

\noindent
{\bf Definition 5.2.}
Let $k_1,k_2,k_3$ satisfy $k_1\geq k_2\geq k_3$.
Take any element $a_1\ot a_2\ot b$ in
$B_{\La_{k_1}}\ot B_{\La_{k_2}}\ot B_{\La_{k_3}}$.
Let
$$
\normalbaselines\mmmth\offinterlineskip
a_1=
	\vcenter{
   \hbox{$\Flect(0.7cm,0.7cm,\hbox{$i_1$})$
         }\vskip-0.4pt
   \hbox{$\Flect(1.4cm,0.7cm,\hbox{$\vdots$})$
         }\vskip-0.4pt
   \hbox{$\Flect(0.7cm,0.7cm,\hbox{$i_{k_1}$})$
         }
                 },
\quad
a_2=
	\vcenter{
   \hbox{$\Flect(0.7cm,0.7cm,\hbox{$i_1^\prime$})$
         }\vskip-0.4pt
   \hbox{$\Flect(1.4cm,0.7cm,\hbox{$\vdots$})$
         }\vskip-0.4pt
   \hbox{$\Flect(0.7cm,0.7cm,\hbox{$i_{k_2}^\prime$})$
         }
                 }
$$
and specify the order on $i_1,\cdots,i_{k_1}$, that is,
fix a bijection $J$
$$
J:\{1,\cdots,k_1\}\ra\{i_1,\cdots,i_{k_1}\}.
$$
We simply call $J$ the order of $a_1$.
Then we define the local index $\ind_{(a_1,J)}(a_2)$
of $a_2$ with respect to $J$
as follows :
\item{1.} Let us write $a_1\ot a_2$ as a two column diagram with dotted boxes
as in the previous sections.
Namely, the first column is $a_1$ and the second column is
$a_2$ (see example below).
Begin from the $J(1)$-th dot in $a_1$.
Look for the dot in $a_2$ which is no lower than this dot.
If there is no such dot in $a_2$, look for the first dot
in $a_2$ from the bottom.
Connect these two dots by a line.
\item{2.} Ignore the dots already connected by a line
and continue the
process 1 for $J(2),\cdots$, $J(k_2)$-th dots in $a_1$.
\item{3.} Define
$\ind_{(a_1,J)}(a_2)$ as the number of lines which connect
dots in $a_1$ with lower dots in $a_2$.
\vskip5mm

\noindent
Next the index $\ind_{(a_1,J)}(b)$ is defined as follows.
By assigning the number $j$ to the dot in $a_2$ which is joined
by a line with the $J(j)$-th dot in $a_1$, we can define the order
of $a_2$.
Let us represent this order by $\bar{J}$,
which gives the bijection
$$
\bar{J}:\{1,\cdots,k_2\}\ra\{i_1^\prime,\cdots,i_{k_2}^\prime\}.
$$
Then we define
$$
\ind_{(a_1,J)}(b)=\ind_{(a_1,J)}(a_2)+\ind_{(a_2,\bar{J})}(b).
$$
\vskip5mm

Let us give an example.
\vskip5mm

\noindent
{\bf Example 5.3.} \
For $sl_4$, consider the element $a_1\ot a_2\ot b$ in
$B_{\La_3}\ot B_{\La_2}\ot B_{\La_1}$ given by
$$
	\normalbaselines\mmmth\offinterlineskip
a_1=	\vcenter{
   \hbox{$\Flect(0.7cm,0.7cm,\hbox{$1$})$
         }\vskip-0.4pt
   \hbox{$\Flect(0.7cm,0.7cm,\hbox{$2$})$
         }\vskip-0.4pt
   \hbox{$\Flect(0.7cm,0.7cm,\hbox{$3$})$
         }
                },
\
a_2=	\vcenter{
   \hbox{$\Flect(0.7cm,0.7cm,\hbox{$1$})$
         }\vskip-0.4pt
   \hbox{$\Flect(0.7cm,0.7cm,\hbox{$3$})$
         }
                },
\
a_2=	\vcenter{
   \hbox{$\Flect(0.7cm,0.7cm,\hbox{$2$})$
         }
                }.
$$
Let us consider the order $J$ given by
$$
J(1)=2,\quad J(2)=3,\quad J(3)=1.
$$
Then the diagram in this case is given by {\bf Figure 5.1}.
{}From this we have
$$
\ind_{(a_1,J)}(a_2)=0,
\quad
\ind_{(a_2,\bar{J})}(b)=\ind_{(a_1,J)}(b)=1.
$$
\vskip5mm

Let us introduce terminology associated with the above definition
of the local index.
\vskip5mm

\noindent
{\bf Definition 5.4.}
We call the line which is drawn in the definition of the local index
the LS-line. A LS-line is called a "down line" if it connects
two dots such that the position of the dot in the left column
is higher than that of the dot in the right column.
Other LS-lines are called "up lines".
We draw a down line in such a way that it first goes right
and then goes up to arrive
at the top horizontal line. Then,
rounding outside the diagram to arrive at the bottom horizontal line,
it finally goes up again to arrive at the ending point (see {\bf Figure 5.1}).
Thus a down line is "winding".
The up line should be similarly drawn.
\vskip1cm

\noindent
{\bf 5.5. Reduction of the statement}\par
Let $b_1\ot\cdots\ot b_l\in B_{\La_{\mu_1}}\ot\cdots\ot B_{\La_{\mu_l}}$
be an element from ${\cal T}(\la,\mu)$. Note that
$$
\normalbaselines\mmmth\offinterlineskip
b_1=
	\vcenter{
   \hbox{$\Flect(0.7cm,0.7cm,\hbox{$1$})$
         }\vskip-0.4pt
   \hbox{$\Flect(1.4cm,0.7cm,\hbox{$\vdots$})$
         }\vskip-0.4pt
   \hbox{$\Flect(0.7cm,0.7cm,\hbox{$\mu_1$})$
         }
                 }
$$
since $b_1\ot\cdots\ot b_l$ is a highest weight element( for $sl_n$).

We define the order $J_1$ of $b_1$ from the bottom up, that is
$$
J_1(i)=\mu_1+1-i.
$$
The order $J_1$ of $b_1$ induces the order $J_2,\cdots,J_l$
of $b_2,\cdots, b_l$, respectively, as in the definition
of the local index. By definition of the local index
and the index of Lascoux-Schutzenberger, we have
$$
\eqalign{
\ind(b_1):=\ind(1)&=0
\cr
\ind(b_i):=\ind(i)&=\ind_{(b_1,J_1)}(b_2)+\cdots+
\ind_{(b_{i-1},J_{i-1})}(b_i).
\cr
}
$$

In order to prove Theorem 4.1, it is sufficient to prove the equation
$$
\ind(b_i^{(j+1)})-\ind(b_i^{(j)})=
-H_{\La_{\mu_j}\La_{\mu_i}}(b_j,b_i^{(j+1)}).
\eqno(5.1)
$$
In fact, summing up (5.1) in $j$ we have
$$
\ind(b_i)-\ind(b_i^{(2)})=
-\sum_{j=2}^{i-1}H_{\La_{\mu_j}\La_{\mu_i}}(b_j,b_i^{(j+1)}).
\eqno(5.2)
$$
Using this and the following lemma we have Theorem 4.1.
\vskip5mm

\noindent
{\bf Lemma 5.6.}
Let us assume $k\geq k^\prime$.
If $a$ is the highest weight element of $B_{\La_k}$,
we have, for any $b\in B_{\La_{k^\prime}}$,
$$
\ind_{(a,J_a)}(b)=-H_{\La_k\La_{k^\prime}}(a,b),
$$
where $J_a$ is the order of $a$ such that $J_a(i)=k+1-i$.
\vskip3mm
\noindent
Proof. Since
$$
\normalbaselines\mmmth\offinterlineskip
a=
	\vcenter{
   \hbox{$\Flect(0.7cm,0.7cm,\hbox{$1$})$
         }\vskip-0.4pt
   \hbox{$\Flect(1.4cm,0.7cm,\hbox{$\vdots$})$
         }\vskip-0.4pt
   \hbox{$\Flect(0.7cm,0.7cm,\hbox{$k$})$
         }
                 }
$$
the lemma follows immediately from the rule of the energy function.
q.e.d.
\vskip5mm

Now to prove (5.1) it is sufficient to prove
\vskip5mm

\noindent
{\bf Proposition 5.7.}
For any order $J$ of $b_{j-1}$ we have
$$
\ind_{(b_{j-1},J)}(b_i^{(j+1)})
-\ind_{(b_{j-1},J)}(b_i^{(j)})
=-H_{\La_{\mu_j}\La_{\mu_i}}(b_j,b_i^{(j+1)}).
$$
\vskip5mm

Finally, in order prove this proposition, it is sufficient to
prove
\vskip5mm

\noindent
{\bf Proposition 5.8.}
Let us assume $k''\geq k'\geq k$.
Taking any element $a_1\ot a_2\ot b$ in
$B_{\La_{k''}}\ot B_{\La_{k'}}\ot B_{\La_{k}}$, we define
$b^\prime$ and $a_2^\prime$ by the isomorphism
$$
\matrix{
B_{\La_{k''}}\ot B_{\La_{k'}}\ot B_{\La_k}&
\simeq&
B_{\La_{k''}}\ot B_{\La_{k}}\ot B_{\La_{k'}}
\cr
a_1\ot a_2\ot b&
\mapsto&
a_1\ot b^\prime\ot a_2^\prime\cr
}.
$$
Then for any order $J$ of $a_1$ we have
$$
\ind_{(a_1,J)}(b)-\ind_{(a_1,J)}(b^\prime)
=-H_{\La_{k^\prime}\La_k}(a_2,b).
$$
\vskip1cm

\noindent
{\bf 5.9. Proof of Proposition 5.8}\par
\vskip5mm

\noindent
{\bf Lemma 5.10.}
For any order $J$ of $a_1$ we have
$$
\ind_{(a_1,J)}(b)-\ind_{(a_1,J)}(a_2)
\geq-H_{\La_{k^\prime}\La_k}(a_2,b).
$$
\vskip3mm
\noindent
Proof.
Note that the left hand side of the above inequality
is nothing but $\ind_{(a_2,\bar{J})}(b)$,
by Definition 5.2.
 Thus we consider
$a_2\ot b$ below.
Let us set
$$
k_2=-H_{\La_{k^\prime}\La_k}(a_2,b).
$$
Then
$$
k_2=\sharp(\hbox{winding H-lines}).
$$
By the rule of drawing the LS-line and H-line, we easily have
(see {\bf Figure 5.2})
$$
\sharp(\hbox{up line})\leq
\sharp(	\hbox{non-winding H-lines})=k-k_2,
$$
which is equivalent to
$$
\sharp(\hbox{down line})=k-\sharp(\hbox{up line})\geq k_2.
$$
q.e.d.
\vskip5mm

By Lemma 5.10 we can write
$$
\ind_{(a_1,J)}(b)-\ind_{(a_1,J)}(a_2)=k_2+l
$$
for some non-negative integer $l$. Then
the proof of Proposition 5.8 reduces to proving the equation
$$
\ind_{(a_1,J)}(b^\prime)-\ind_{(a_1,J)}(a_2)=l.
\eqno(5.3)
$$
The proof of (5.3) is divided into several steps.
\vskip5mm

\noindent
{Step 1.}
\par
\noindent
Let us consider the following situation ({\bf Figure 5.3}).
\item{1.} In {\bf Figure 5.3}, the line starting from
$
A\hskip1mm\vcenter{
   \hbox{$\Flect(0.3cm,0.3cm,\hbox{$\bullet$})$
         }
                 }
$
is a down line and below it there are no
$
\vcenter{
   \hbox{$\Flect(0.3cm,0.3cm,\hbox{$\bullet$})$
         }
                 }
$
from which a down line starts.
\item{2.} The line entering
$
\vcenter{
   \hbox{$\Flect(0.3cm,0.3cm,\hbox{$\bullet$})$
         }
                 }B
$
is a down line and above it there are no
$
\vcenter{
   \hbox{$\Flect(0.3cm,0.3cm,\hbox{$\bullet$})$
         }
                 }
$
into which a down line enters.
\vskip5mm

\noindent
{\bf Lemma 5.11.}
$
\vcenter{
   \hbox{$\Flect(0.3cm,0.3cm,\hbox{$\bullet$})$
         }
                 }B
$
is below
$
A\hskip1mm\vcenter{
   \hbox{$\Flect(0.3cm,0.3cm,\hbox{$\bullet$})$
         }
                 }.
$
\vskip3mm
\noindent
Proof.
\item{i.} If the LS-line starting from
$
A\hskip1mm\vcenter{
   \hbox{$\Flect(0.3cm,0.3cm,\hbox{$\bullet$})$
         }
                 }
$
and that ending in
$
\vcenter{
   \hbox{$\Flect(0.3cm,0.3cm,\hbox{$\bullet$})$
         }
                 }B
$
are the same, the lemma is obvious by the drawing rule
of LS-lines.
\item{ii.} Suppose that
the LS-line starting from
$
A\hskip1mm\vcenter{
   \hbox{$\Flect(0.3cm,0.3cm,\hbox{$\bullet$})$
         }
                 }
$
and the LS-line ending in
$
\vcenter{
   \hbox{$\Flect(0.3cm,0.3cm,\hbox{$\bullet$})$
         }
                 }B
$
are distinct.
By the definition of
$
A\hskip1mm\vcenter{
   \hbox{$\Flect(0.3cm,0.3cm,\hbox{$\bullet$})$
         }
                 }
$
the starting point of the LS-line ending in
$
\vcenter{
   \hbox{$\Flect(0.3cm,0.3cm,\hbox{$\bullet$})$
         }
                 }B
$
is above
$
A\hskip1mm\vcenter{
   \hbox{$\Flect(0.3cm,0.3cm,\hbox{$\bullet$})$
         }
                 }.
$
Now suppose that
$
\vcenter{
   \hbox{$\Flect(0.3cm,0.3cm,\hbox{$\bullet$})$
         }
                 }B
$
is above
$
A\hskip1mm\vcenter{
   \hbox{$\Flect(0.3cm,0.3cm,\hbox{$\bullet$})$
         }
                 }.
$
Then the order of the line connected to
$
\vcenter{
   \hbox{$\Flect(0.3cm,0.3cm,\hbox{$\bullet$})$
         }
                 }B
$
is smaller than that of
$
A\hskip1mm\vcenter{
   \hbox{$\Flect(0.3cm,0.3cm,\hbox{$\bullet$})$
         }
                 }.
$
In fact, if this were not the case, the LS-line starting from
$
A\hskip1mm\vcenter{
   \hbox{$\Flect(0.3cm,0.3cm,\hbox{$\bullet$})$
         }
                 }
$
would stop at
$
\vcenter{
   \hbox{$\Flect(0.3cm,0.3cm,\hbox{$\bullet$})$
         }
                 }B
$
or below it (see {\bf Figure 5.4}).
This contradicts the assumption that the line from
$
A\hskip1mm\vcenter{
   \hbox{$\Flect(0.3cm,0.3cm,\hbox{$\bullet$})$
         }
                 }
$
is down.
Then there is no place where the line starting from
$
A\hskip1mm\vcenter{
   \hbox{$\Flect(0.3cm,0.3cm,\hbox{$\bullet$})$
         }
                 }
$
ends, since this line must end between the bottom and
$
\vcenter{
   \hbox{$\Flect(0.3cm,0.3cm,\hbox{$\bullet$})$
         }
                 }B
$
by the definition of
$
\vcenter{
   \hbox{$\Flect(0.3cm,0.3cm,\hbox{$\bullet$})$
         }
                 }B.
$
q.e.d.
\vskip5mm

\noindent
{\bf Lemma 5.12.}
Between
$
A\hskip1mm\vcenter{
   \hbox{$\Flect(0.3cm,0.3cm,\hbox{$\bullet$})$
         }
                 }
$
and
$
\vcenter{
   \hbox{$\Flect(0.3cm,0.3cm,\hbox{$\bullet$})$
         }
                 }B
$
there exists
$
C\hskip1mm\vcenter{
   \hbox{$\Flect(0.3cm,0.3cm,\hbox{$$})$
         }
                 }
$
in $a_2$ whose bottom line does not intersect any LS-line
(see {\bf Figure 5.5}).
$
C\hskip1mm\vcenter{
   \hbox{$\Flect(0.3cm,0.3cm,\hbox{$$})$
         }
                 }
$
may coincide with
$
A\hskip1mm\vcenter{
   \hbox{$\Flect(0.3cm,0.3cm,\hbox{$$})$
         }
                 }
$,
but must be distinct from
$
\vcenter{
   \hbox{$\Flect(0.3cm,0.3cm,\hbox{$$})$
         }
                 }B.
$
\vskip3mm
\noindent
Proof.
Suppose that such a
$
C\hskip1mm\vcenter{
   \hbox{$\Flect(0.3cm,0.3cm,\hbox{$$})$
         }
                 }
$
does not exist. In this case the situation is like that in
{\bf Figure 5.6}. We proceed with the proof following this picture.
The order of the LS-line $\alpha_1-\alpha_5$ is as follows:
\item{$\bullet$} The order of $\alpha_3$ is smaller than that of $\alpha_1$.
\item{$\bullet$} The order of $\alpha_4$ is smaller than that of $\alpha_3$.
\item{$\bullet$} The order of $\alpha_5$ is smaller than that of $\alpha_4$.
\item{$\bullet$} The order of $\alpha_2$ is smaller than that of $\alpha_5$.
\par
\noindent
Hence
the order of $\alpha_2$ is smaller than that of $\alpha_1$.
This contradicts the fact that the end point of the line $\alpha_1$
is lower than the position of
$
\vcenter{
   \hbox{$\Flect(0.3cm,0.3cm,\hbox{$\bullet$})$
         }
                 }B.
$
Note that the general case is to consider the $m$ LS-lines
$\alpha_1-\alpha_m$ and is proved analogously to the case above.
q.e.d.
\vskip5mm

Now using Lemma 5.12 we modify our two column diagram in the following manner.
Take any
$
C\hskip1mm\vcenter{
   \hbox{$\Flect(0.3cm,0.3cm,\hbox{$$})$
         }
                 }
$
as in Lemma 5.12.
We cut the diagram by the bottom horizontal line of
$
C\hskip1mm\vcenter{
   \hbox{$\Flect(0.3cm,0.3cm,\hbox{$$})$
         }
                 }
$
and put the lower half of this diagram over
the upper half (see {\bf Figure 5.7}).
We can naturally transplant the LS-lines to this new diagram
(see {\bf Figure 5.8}). We also call those new lines LS-lines.
Then the LS-lines for this new diagram have the following
properties:
\item{$\bullet$} The new LS-lines go up only.
\item{$\bullet$} Down lines in the old diagram correspond bijectively
to new lines which intersect the top horizontal line of the
old diagram.
\par
\noindent
In the following we consider only the new diagram, unless otherwise
stated, and call the
top horizontal line of the old diagram the "{\bf bold line}" in the new
diagram.
\vskip1cm

\noindent
{Step 2.}
\par
\noindent
To each
$
\vcenter{
   \hbox{$\Flect(0.3cm,0.3cm,\hbox{$\bullet$})$
         }
                 }
$
in the right column and
part of
$
\vcenter{
   \hbox{$\Flect(0.3cm,0.3cm,\hbox{$\bullet$})$
         }
                 }
$
in the left column
 we attach the symbol $u$ or $d$ by the following rule:
\item{1.} Begin with the lowest
$
\vcenter{
   \hbox{$\Flect(0.3cm,0.3cm,\hbox{$\bullet$})$
         }
                 }
$
in the right column.
In the left column look for the nearest
$
\vcenter{
   \hbox{$\Flect(0.3cm,0.3cm,\hbox{$\bullet$})$
         }
                 }
$
which is no higher than this
$
\vcenter{
   \hbox{$\Flect(0.3cm,0.3cm,\hbox{$\bullet$})$
         }
                 }
$,
and connect these two
$
\vcenter{
   \hbox{$\Flect(0.3cm,0.3cm,\hbox{$\bullet$})$
         }
                 }
$
with a line. If this line crosses the bold line then we attach
the symbol $d$ to the original
$
\vcenter{
   \hbox{$\Flect(0.3cm,0.3cm,\hbox{$\bullet$})$
         }
                 }
$
in the right column as :
$
\vcenter{
   \hbox{$\Flect(0.3cm,0.3cm,\hbox{$\bullet$})$
         }
                 }d.
$
Otherwise we attach the symbol $u$ as :
$
\vcenter{
   \hbox{$\Flect(0.3cm,0.3cm,\hbox{$\bullet$})$
         }
                 }u.
$
\item{2.} Ignore the pair already joined by a line
and continue the process 1 for the 2nd,3rd,...
$
\vcenter{
   \hbox{$\Flect(0.3cm,0.3cm,\hbox{$\bullet$})$
         }
                 }
$
from the bottom successively.
\item{3.}
For
$
\vcenter{
   \hbox{$\Flect(0.3cm,0.3cm,\hbox{$\bullet$})$
         }
                 }
$
in the left column which are connected by a line to
some
$
\vcenter{
   \hbox{$\Flect(0.3cm,0.3cm,\hbox{$\bullet$})$
         }
                 }
$
in the right column, we attach $u$ or $d$ to match the
corresponding symbol of
$
\vcenter{
   \hbox{$\Flect(0.3cm,0.3cm,\hbox{$\bullet$})$
         }
                 }
$
in the right column.
\item{4.} Finally we attach $v$ to any
$
\vcenter{
   \hbox{$\Flect(0.3cm,0.3cm,\hbox{$\bullet$})$
         }
                 }
$
which has neither a $u$ nor a $d$ as :
$
v\hskip1mm\vcenter{
   \hbox{$\Flect(0.3cm,0.3cm,\hbox{$\bullet$})$
         }
                 }.
$
\vskip3mm

\noindent
In processes $1$ and $2$ above, the lines never "wind".
In fact let us consider process $1$.
There exists an up LS-line whose end point is the lowest
$
\vcenter{
   \hbox{$\Flect(0.3cm,0.3cm,\hbox{$\bullet$})$
         }
                 }
$
in the right column.
Let us name its starting point $A$. Then the chosen
$
\vcenter{
   \hbox{$\Flect(0.3cm,0.3cm,\hbox{$\bullet$})$
         }
                 }
$
in the left column in process $1$ is not below $A$.
Hence the line drawn in process $1$ cannot wind.
The same reasoning is applied to each case of process $2$.
\vskip3mm

\noindent
{\bf Definition 5.13.}
We refer to the line drawn above the $H$-line for the new diagram.
\vskip3mm

\noindent
After finishing this process, the symbols $u$ and $d$ have the
following properties (see {\bf Figure 5.9}):
\item{1.} All
$
\vcenter{
   \hbox{$\Flect(0.3cm,0.3cm,\hbox{$\bullet$})$
         }
                 }d
$
in the right column are over the bold line.
\item{2.}All
$
\vcenter{
   \hbox{$\Flect(0.3cm,0.3cm,\hbox{$\bullet$})$
         }
                 }d
$
in the left column are under the bold line.

\item{3.} There is no
$
v\hskip1mm\vcenter{
   \hbox{$\Flect(0.3cm,0.3cm,\hbox{$\bullet$})$
         }
                 }
$
which is under the bold line and is over the lowest
$
d\hskip1mm\vcenter{
   \hbox{$\Flect(0.3cm,0.3cm,\hbox{$\bullet$})$
         }
                 }.
$
\item{4.} Let
$
u\hskip1mm\vcenter{
   \hbox{$\Flect(0.3cm,0.3cm,\hbox{$\bullet$})$
         }
                 }=A
$
be the $j$-th
$
u\hskip1mm\vcenter{
   \hbox{$\Flect(0.3cm,0.3cm,\hbox{$\bullet$})$
         }
                 }
$
counting from the bottom (we count only
$
u\hskip1mm\vcenter{
   \hbox{$\Flect(0.3cm,0.3cm,\hbox{$\bullet$})$
         }
                 }'s
$).
If $A$ sits lower than the lowest
$
d\hskip1mm\vcenter{
   \hbox{$\Flect(0.3cm,0.3cm,\hbox{$\bullet$})$
         }
                 },
$
then the number of LS-lines starting from
$
\vcenter{
   \hbox{$\Flect(0.3cm,0.3cm,\hbox{$$})$
         }
                 }
$
in the left column whose positions are no higher than $A$
is greater than or equal to $j$.
\par
\vskip2mm

\noindent
Properties 1-3 are obvious.
Let us prove property 4.
Let us name
$
\vcenter{
   \hbox{$\Flect(0.3cm,0.3cm,\hbox{$\bullet$})$
         }
                 }
$
from bottom to up as $A_1,A_2,\cdots$.
Let $A_N$ be the highest $A_i$ which is lower than the lowest
$
d\hskip1mm\vcenter{
   \hbox{$\Flect(0.3cm,0.3cm,\hbox{$\bullet$})$
         }
                 }.
$
Let us also name
$
\vcenter{
   \hbox{$\Flect(0.3cm,0.3cm,\hbox{$\bullet$})$
         }
                 }u
$
from the bottom up as $B_1,B_2,\cdots$.
Let us first prove that property 4 holds for $j=N$.
Note that the position of $B_N$ is not lower than that of $A_N$
and is lower than the lowest
$
d\hskip1mm\vcenter{
   \hbox{$\Flect(0.3cm,0.3cm,\hbox{$\bullet$})$
         }
                 }.
$
In fact if $B_N$ is lower than $A_N$, then there are already
$N$
$
u\hskip1mm\vcenter{
   \hbox{$\Flect(0.3cm,0.3cm,\hbox{$\bullet$})$
         }
                 }
$'s
below $A_N$. This is impossible.
By the drawing order of $H$ lines $B_N$, should be below the lowest
$
d\hskip1mm\vcenter{
   \hbox{$\Flect(0.3cm,0.3cm,\hbox{$\bullet$})$
         }
                 }.
$
Moreover, there are no
$
v\hskip1mm\vcenter{
   \hbox{$\Flect(0.3cm,0.3cm,\hbox{$\bullet$})$
         }
                 }
$
which are above $A_N$ and not above $B_N$
by the drawing rule of $H$ lines.
Therefore no starting point of $H$ lines which end at
$B_1,B_2,\cdots$ can be above $A_N$.
Thus property 4 holds for $j=N$.

Now suppose that property 4 holds for all
$j^\prime(j\leq j^\prime\leq N)$
and does not hold for $j-1$.
Then there should be a
$
v\hskip1mm\vcenter{
   \hbox{$\Flect(0.3cm,0.3cm,\hbox{$\bullet$})$
         }
                 }
$
between $A_{j-1}$ and $A_j$ from which an $H$ line starts.
Let $C$ be the lowest
$
v\hskip1mm\vcenter{
   \hbox{$\Flect(0.3cm,0.3cm,\hbox{$\bullet$})$
         }
                 }
$
which sits between $A_{j-1}$ and $A_j$.
Then $B_{j-1}$ is no lower than $A_{j-1}$ but is lower than $C$.
In fact the position of $B_{j-1}$ cannot be below that of
$A_{j-1}$ by the same reasoning as in the case of $A_N$ and $B_N$.
By the drawing rule of $H$ lines (or $u$), $B_{j-1}$ must be below $C$.
As a consequence all starting points of $H$ lines ending at
$B_1,\cdots,B_{j-1}$ are under $A_{j-1}$.
Thus property 4 holds for $j-1$. This is a contradiction. q.e.d.
\vskip5mm

\noindent
Note that an old $H$ line drawn by the above order
is naturally transplanted as a new $H$ line.
By this transplant the winding $H$ line is interpreted
as the new $H$ line which crosses the "bold line".
\vskip5mm
\noindent
{Step 3.}
\par
\noindent
Let us set
$$
\eqalign{
&\hbox{
The number of
$
\vcenter{
   \hbox{$\Flect(0.3cm,0.3cm,\hbox{$\bullet$})$
         }
                 }u
$
over the bold line in the right column}
=N_R^{\hbox{up}}(u).
\cr
&\hbox{
The number of
$
\vcenter{
   \hbox{$\Flect(0.3cm,0.3cm,\hbox{$\bullet$})$
         }
                 }u
$
under the bold line in the right column}
=N_R^{\hbox{down}}(u).
}
$$
Then we have the following statements.
$$
\eqalign{
&\hbox{
The number of
$
u\hskip1mm\vcenter{
   \hbox{$\Flect(0.3cm,0.3cm,\hbox{$\bullet$})$
         }
                 }
$
over the bold line in the left column}
=N_R^{\hbox{up}}(u).
\cr
&\hbox{
The number of
$
u\hskip1mm\vcenter{
   \hbox{$\Flect(0.3cm,0.3cm,\hbox{$\bullet$})$
         }
                 }
$
under the bold line in the left column}
=N_R^{\hbox{down}}(u).
\cr
&\hbox{
The number of
$
d\hskip1mm\vcenter{
   \hbox{$\Flect(0.3cm,0.3cm,\hbox{$\bullet$})$
         }
                 }
$
under the bold line in the left column}
=k_2.
}
$$
Here $k_2$ is the number of winding $H$ lines in $a_2\ot b$
in our assumption.
Hence we have
$$
\sharp(
u\hskip1mm\vcenter{
   \hbox{$\Flect(0.3cm,0.3cm,\hbox{$\bullet$})$
         }
                 }
\hbox{ under the bold line})
+
\sharp(
d\hskip1mm\vcenter{
   \hbox{$\Flect(0.3cm,0.3cm,\hbox{$\bullet$})$
         }
                 }
\hbox{ under the bold line})
=N_R^{\hbox{down}}(u)+k_2.
\eqno{(5.4)}
$$
On the other hand,
\item{$\bullet$}
The number of LS-lines whose starting points are  under the bold line
and end points are over the bold line is $k_2+l$ by assumption.
\item{$\bullet$}The number of LS-lines whose starting and end points
are both under the bold line is $N_R^{\hbox{down}}(u)$.
\par
\noindent
Consequently,
$$\eqalign{
&\hbox{(the number of LS-lines whose starting points are under the bold
line)} \cr
&= k_2+l+N_R^{\hbox{down}}(u).
}\eqno(5.5)
$$
\vskip2mm

To this point we have considered LS-lines starting
from the left column $a_2$. Now we consider the LS-lines
coming into $a_2$ from $a_1$.
\vskip3mm
\noindent
{\bf Lemma 5.14.}
If we remove all
$
v\hskip1mm\vcenter{
   \hbox{$\Flect(0.3cm,0.3cm,\hbox{$\bullet$})$
         }
                 }
$
from the left column, among the LS-lines going into the left column,
$l$ LS-lines "newly cross the bold line".
Here, for example, in
{\bf Figure 5.10}, the line $\alpha_3$  "newly crosses the bold line"
after removing
$
v\hskip1mm\vcenter{
   \hbox{$\Flect(0.3cm,0.3cm,\hbox{$\bullet$})$
         }
                 }.
$
\vskip3mm
\noindent
Proof. We first note the following property on the movement
of the end position of LS-lines while removing
$
v\hskip1mm\vcenter{
   \hbox{$\Flect(0.3cm,0.3cm,\hbox{$\bullet$})$
         }
                 }'s
$:
$$
\eqalign{
&\hbox{If we remove
$
v\hskip1mm\vcenter{
   \hbox{$\Flect(0.3cm,0.3cm,\hbox{$\bullet$})$
         }
                 }
$,
then the end positions of LS-lines in the left column
shift upward only}
\cr
&\hbox{under the periodic boundary condition.}
}
\eqno(5.6)
$$
(see {\bf Figure 5.11}).
\noindent
More strongly, we can prove
\vskip3mm
\noindent
{\bf Lemma 5.15.} There does not occur any winding
of lines (in the new diagram) after removing
$
v\hskip1mm\vcenter{
   \hbox{$\Flect(0.3cm,0.3cm,\hbox{$\bullet$})$
         }
                 }
$
compared with the situation before removing
$
v\hskip1mm\vcenter{
   \hbox{$\Flect(0.3cm,0.3cm,\hbox{$\bullet$})$
         }
                 }.
$
\vskip2mm
\noindent
Proof.
In this proof we consider $a_2\ot b$. So the starting point
of an LS-line is again a point in $a_2$.
Let us consider the highest
$
\hskip1mm\vcenter{
   \hbox{$\Flect(0.3cm,0.3cm,\hbox{$\bullet$})$
         }
                 }d
$
(say $D$) over the bold line and $L$ the (horizontal) overline
of the box $D$ ({\bf Figure 5.12}).
Then any
$
\hskip1mm\vcenter{
   \hbox{$\Flect(0.3cm,0.3cm,\hbox{$\bullet$})$
         }
                 }
$
in the left column which is above the bold line and under $L$
is assigned the symbol $u$ by definition of the symbol
$
\hskip1mm\vcenter{
   \hbox{$\Flect(0.3cm,0.3cm,\hbox{$\bullet$})$
         }
                 }d.
$
Moreover any
$
\hskip1mm\vcenter{
   \hbox{$\Flect(0.3cm,0.3cm,\hbox{$\bullet$})$
         }
                 }
$
above $L$ in the right column is assigned the symbol $u$
by definition of $L$.
Let $A_1,\cdots,A_N$ be the starting points of LS-lines over $L$
numbered from the bottom up and
$l_1,\cdots,l_n$ the corresponding LS-lines.
Let $B_1,\cdots,B_N$ be the end points of $\{l_i\}$
numbered from the bottom up ({\bf Figure 5.13}).

In the following we use the notation $A\leq B$ for two boxes $A$ and $B$
meaning the position of $A$ is equal to or lower than that of $B$.
The symbol $A<B$ is similarly defined.

Note that $A_j\leq B_j$. In fact if $A_j>B_j$, then there are
$j$ LS-lines going from $\{A_i\vert A_i<B_j<A_j\}$ to
$B_1,\cdots,B_j$. This is impossible.
Note also the following.
The
$
\hskip1mm\vcenter{
   \hbox{$\Flect(0.3cm,0.3cm,\hbox{$\bullet$})$
         }
                 }
$
in the left column which is joined by an $H$ line to some
$
\hskip1mm\vcenter{
   \hbox{$\Flect(0.3cm,0.3cm,\hbox{$\bullet$})$
         }
                 }u
$
over $L$ is always over $L$ by definition of the rule of drawing
$H$ lines in the present order.
So we can consider the part over $L$ separately.
Temporarilly, we consider only this part unless otherwise stated.
Then, by the fiirst part of the proof of Lemma 3.16,
the set of
$
u\hskip1mm\vcenter{
   \hbox{$\Flect(0.3cm,0.3cm,\hbox{$\bullet$})$
         }
                 }
$
(over $L$) is independent of the order of drawing $H$ lines
(whose starting points are over $L$).
Now let us prove the following claim
$$
\eqalign{
&\hbox{There exists a subset $\{C_1,\cdots,C_N\}$ of
$\{
u\hskip1mm\vcenter{
   \hbox{$\Flect(0.3cm,0.3cm,\hbox{$\bullet$})$
         }
                 }
\}
$
such that $C_1<\cdots<C_N$ and}
\cr
&\hbox{$A_1\leq C_1,\cdots,A_N\leq C_N$}.
}
\eqno(5.7)
$$
In fact let us define the order $1,2,\cdots,N$ to $B_N,\cdots,B_1$.
The order $N+1,N+2,\cdots$ is irrelevant to the following argument
and arbitrarily assigned.
Let $\bar{B}_N,\cdots,\bar{B}_1$ be the partners of
$B_N,\cdots,B_1$ by $H$ lines in the presently specified order.
Then by definition
$$
\bar{B}_1<\cdots<\bar{B}_N.
$$
Since $A_j\leq B_j$, as we already proved, $A_j\leq \bar{B}_j$.
Hence we can choose $C_j=\bar{B}_j$.
Thus (5.7) is proved.

Let us return to the proof of Lemma 5.15.
{}From this point, we again consider the entire diagram.
Under the bold line in the left column we have
$$
\eqalign{
&\sharp(
u\hskip1mm\vcenter{
   \hbox{$\Flect(0.3cm,0.3cm,\hbox{$\bullet$})$
         }
                 }
\hbox{ and }
d\hskip1mm\vcenter{
   \hbox{$\Flect(0.3cm,0.3cm,\hbox{$\bullet$})$
         }
                 }
)
-\sharp(\hbox{starting points of LS-lines})
\cr
&
=(k_2+N_R^{\hbox{down}}(u))-(k_2+l+N_R^{\hbox{down}}(u))
\cr
&=-l.
}
$$
Hence, over the bold line in the left column,
$$
\sharp(
u\hskip1mm\vcenter{
   \hbox{$\Flect(0.3cm,0.3cm,\hbox{$\bullet$})$
         }
                 }
\hbox{ and }
d\hskip1mm\vcenter{
   \hbox{$\Flect(0.3cm,0.3cm,\hbox{$\bullet$})$
         }
                 }
)
-\sharp(\hbox{starting points of LS-lines})=l.
$$
Since, in the left column as a whole, the following equation must hold
$$
\sharp(
u\hskip1mm\vcenter{
   \hbox{$\Flect(0.3cm,0.3cm,\hbox{$\bullet$})$
         }
                 }
\hbox{ and }
d\hskip1mm\vcenter{
   \hbox{$\Flect(0.3cm,0.3cm,\hbox{$\bullet$})$
         }
                 }
)
=\sharp(\hbox{starting points of LS-lines}).
$$
Let us prove the impossibility of winding by considering two cases.
\vskip1mm
\noindent
Case 1. The case that the line (say $L$) which winds for the first time
(in the order of drawing LS-lines) "newly crosses the bold line".
\vskip2mm
\noindent
Suppose that this occurs.
Then every
$
u\hskip1mm\vcenter{
   \hbox{$\Flect(0.3cm,0.3cm,\hbox{$\bullet$})$
         }
                 }
$
and
$
d\hskip1mm\vcenter{
   \hbox{$\Flect(0.3cm,0.3cm,\hbox{$\bullet$})$
         }
                 }
$
over the bold line is occupied by some LS-lines
before drawing $L$.
Hence $l$ lines already have "newly crossed the bold line"
before drawing $L$, and $L$ is the $l+1$-th line which
"newly crosses the bold line".
But this is impossible unless a line already has wound
before drawing $L$ by property 4 in Step 2.
\vskip2mm

\noindent
Case 2. The case that the line (say $L$) which winds for the first time
does not "newly cross the bold line".
\vskip2mm
\noindent
In this case $l+1$ lines should "newly cross the bold line"
before drawing $L$.
In fact let $A$ be the starting point of $L$ before removing
$
v\hskip1mm\vcenter{
   \hbox{$\Flect(0.3cm,0.3cm,\hbox{$\bullet$})$
         }
                 }.
$
Then $A$ should be
$
v\hskip1mm\vcenter{
   \hbox{$\Flect(0.3cm,0.3cm,\hbox{$\bullet$})$
         }
                 }
$.
If this is not the case $L$ cannot wind.
If there exists
$
u\hskip1mm\vcenter{
   \hbox{$\Flect(0.3cm,0.3cm,\hbox{$\bullet$})$
         }
                 }
$
or
$
d\hskip1mm\vcenter{
   \hbox{$\Flect(0.3cm,0.3cm,\hbox{$\bullet$})$
         }
                 }
$
below $A$ and above $L$
which is not occupied before drawing $L$,
then $L$ cannot wind by (5.7).
Hence before drawing $L$ every
$
u\hskip1mm\vcenter{
   \hbox{$\Flect(0.3cm,0.3cm,\hbox{$\bullet$})$
         }
                 }
$
and
$
d\hskip1mm\vcenter{
   \hbox{$\Flect(0.3cm,0.3cm,\hbox{$\bullet$})$
         }
                 }
$
over the bold line must become the starting points of some LS-lines.
Hence $l+1$ lines should already "newly cross the bold line"
before drawing $L$. But this is impossible by the same reasoning
as that in Case 1.
q.e.d.
\vskip3mm

We can now finish the proof of Lemma 5.14.
In fact it follows from property 4 in Step 2, (5.6)
and Lemma 5.15 that
the number
of LS-lines "newly crossing the bold line" is
$$
k_2+l+N_R^{\hbox{down}}(u)-(k_2+N_R^{\hbox{down}}(u))=l
$$
after removing
$
v\hskip1mm\vcenter{
   \hbox{$\Flect(0.3cm,0.3cm,\hbox{$\bullet$})$
         }
                 }.
$
q.e.d.
\vskip1cm

\noindent
{Step 4.}
\par
\noindent
Let us prove (5.3).
We consider the old diagram here
(i.e., not that one with the lower and upper halves reversed).
When we remove all
$
v\hskip1mm\vcenter{
   \hbox{$\Flect(0.3cm,0.3cm,\hbox{$\bullet$})$
         }
                 }
$
from $a_2$, the following four changes can occur for LS-lines connecting
$a_1$ and $a_2$ (the following are
$
\vcenter{
   \hbox{$\Flect(0.3cm,0.3cm,\hbox{$$})$
         }
                 }
$
in $a_1$ and the LS-lines coming from it):
\vskip2mm
\noindent
\item{(1)}
$
\vcenter{
   \hbox{$\Flect(0.3cm,0.3cm,\hbox{$$})$
         }
                 }\!\!\!\!\nearrow\hbox{ up line}
$
$\Longrightarrow$
$
\vcenter{
   \hbox{$\Flect(0.3cm,0.3cm,\hbox{$$})$
         }
                 }\!\!\!\!\nearrow\hbox{ up line.}
$
\item{(2)}
$
\vcenter{
   \hbox{$\Flect(0.3cm,0.3cm,\hbox{$$})$
         }
                 }\!\!\!\!\nearrow\hbox{ up line}
$
$\Longrightarrow$
$
\vcenter{
   \hbox{$\Flect(0.3cm,0.3cm,\hbox{$$})$
         }
                 }\!\!\!\!\searrow\hbox{ down line.}
$
\item{(3)}
$
\vcenter{
   \hbox{$\Flect(0.3cm,0.3cm,\hbox{$$})$
         }
                 }\!\!\!\!\searrow\hbox{ down line}
$
$\Longrightarrow$
$
\vcenter{
   \hbox{$\Flect(0.3cm,0.3cm,\hbox{$$})$
         }
                 }\!\!\!\!\nearrow\hbox{ up line.}
$
\item{(4)}
$
\vcenter{
   \hbox{$\Flect(0.3cm,0.3cm,\hbox{$$})$
         }
                 }\!\!\!\!\searrow\hbox{ down line}
$
$\Longrightarrow$
$
\vcenter{
   \hbox{$\Flect(0.3cm,0.3cm,\hbox{$$})$
         }
                 }\!\!\!\!\searrow\hbox{ down line.}
$
\par
\vskip2mm
\noindent
In truth, case $(3)$ cannot occur by the rule  for $H$-lines.
\vskip5mm

\noindent
{\bf Lemma 5.16.}
If and only if the situation is that of case $(2)$,
the LS-line starting from
$
\vcenter{
   \hbox{$\Flect(0.3cm,0.3cm,\hbox{$$})$
         }
                 }
$
in $a_1$ will "newly cross the bold line" after removing all
$
v \hskip1mm \vcenter{
   \hbox{$\Flect(0.3cm,0.3cm,\hbox{$\bullet$})$
         }
                 }
$
from $a_2$.
\vskip2mm
\noindent
Proof. This is obvious from (5.6) and the definition of
the word "newly cross the bold line".
q.e.d.
\vskip5mm

\noindent
Lemma 5.16 states
$$
{\rm ind}_{(a_1,J)}(b')={\rm ind}_{(a_1,J)}(a_2)+
\sharp \{ \hbox{LS-line of type (2)}  \}.
$$
Combining this with Lemma 5.14, we have (5.3).
q.e.d.
\vskip1cm

\vfill\eject

\noindent
{\fontB 6. Summary and Discussion}\par
\vskip2mm
\noindent
In this paper we have given an expression for the charge,
or more precisely the index of
Lascoux-Schutzenberger, in terms of the energy function.
Since the energy function is determined from the crystals
over quantum affine algebra, our theorem provides a
representation theoretical meaning of the charge.
As a corollary of our description of charge, we proved Kirillov's
conjecture on the relation between Kostka polynomials and the branching
functions of the vacuum representation of $\widehat{sl_n}$ w.r.t. $sl_n$.

While preparing the manuscript we came to know about the paper [LLT]
by Lascoux et al. In this paper the charge is also described
in terms of crystal graphs. The relation between their description
and ours is still under investigation.
Here we only remark that their $d_i(t)$ defining the statistics $d(t)$
looks similar to the energy function, in particular in the case
$\mu=(k^{n+1})$ for $sl_{n+1}$. However, they are not identical
since $d_i(t)$
is determined from the $i$ string through $t$, while the energy function
can actually depend on other data. In any case it is an interesting
problem to clarify the relation between the two descriptions of charge.

After this manuscript is completed, there appeared a paper
by Dasmahapatra [D] that containes some results similar to ours.
\vskip5mm

\noindent
{\bf Acknowledgement}\par
\noindent
We wish to thank A.N. Kirillov for stimulating discussions
and informing us of the paper [LLT].

\vfill\eject
%Nakayashiki personal macros
\def\ot{\otimes}
\def\ra{\longrightarrow}
\def\et{\tilde{e}}
\def\ft{\tilde{f}}
\def\La{\Lambda}
\def\la{\lambda}
\def\mapright#1{\smash{\mathop{\longrightarrow}\limits^{#1}}}

\noindent
{\fontB Appendix \ - Fundamentals of Crystal Theory }\par
\noindent
We briefly give here precise definitions of notions and notation
used in the main text. For more details we refer to [KMN1,K3].
The crystals defined below are called seminormal in [K3]
and are the same as the weighted crystal in [KMN1].
We restrict ourselves to the case of $A^{(1)}_n$.

We denote by $\{\alpha_0,\cdots,\alpha_n\}$ the set of simple
roots, by $\{h_0,\cdots,h_n\}$ the set of simple coroots
so that $(<h_i,\alpha_j>)$ is the generalized Cartan matrix
of the $A^{(1)}_n$ type, by $\{\La_0,\cdots,\La_n\}$ the set of
fundamental weights, by $\delta=\alpha_0+\cdots+\alpha_n$ the
null root, by $P={\bf Z}\La_0+\cdots+{\bf Z}\La_n+{\bf Z}\delta$
the weight lattice, and by
$P^\ast={\bf Z}h_0+\cdots+{\bf Z}h_n+{\bf Z}d$ the dual weight lattice.
The pairing of $P$ and $P^\ast$ is given by
$$
\langle h_i,\La_j \rangle=\delta_{ij},\quad
\langle d, \delta \rangle=1,\quad
\hbox{others}=0.
$$
\vskip3mm

\noindent
{\bf Definition A.1.}
A crystal $B$ (for $\widehat{sl_{n+1}}$) is a set with a map
$$
wt:B\ra P,\
\et_i,\ft_i:B\sqcup \{0\}\ra B\sqcup \{0\},\
\varepsilon_i,\varphi_i:B\ra{\bf Z}_{\geq0},\
i=0,\cdots,n
$$
satisfying the following axioms.
\item{1.} If $b\in B$ and $\et_ib\in B$, then $wt(\et_ib)=wtb+\alpha_i$.
\item{2.} If $b\in B$ and $\ft_ib\in B$, then $wt(\ft_ib)=wtb-\alpha_i$.
\item{3.} For $b\in B$,
$\varphi_i(b)=\max\{n\geq0\vert\ft_i^nb\neq0\}$,
$\varepsilon_i(b)=\max\{n\geq0\vert\et_i^nb\neq0\}$.
\item{4.} $\varphi_i(b)-\varepsilon_i(b)=\langle h_i, wtb\rangle$.
\item{5.} $\et_i0=\ft_i0=0$.
\item{6.} For $b_1,b_2\in B$, $\ft_ib_1=b_2$
is equivalent to $b_1=\et_ib_2$.
\par
\vskip3mm

\noindent
Let us set $P_{cl}=P/{\bf Z}\delta$.
A crystal similarly defined for $P_{cl}$ rather than $P$ is called
a classical crystal. In the main text we use the term crystal
to refer to a classical crystal as long as there is
no chance of confusion.
A crystal similarly defined for
$\bar{P}={\bf Z}\La_1+\cdots+{\bf Z}\La_n$ and
$\et_i,\ft_i,\varphi_i,\varepsilon_i(1\leq i\leq n)$
is a crystal for $sl_{n+1}$.
For a subset $J$ of $\{0,1,\cdots,n\}$,
a crystal $B$ (for $\widehat{sl_{n+1}}$)
is called a $J$ crystal if the index of
$\et_i,\ft_i,\varphi_i,\varepsilon_i$ is restricted to $J$.
In particular, for $J=\{1,\cdots,n\}$, the $J$- crystal is sometimes called
a $sl_{n+1}$ crystal.

If $(L,B)$ is a crystal base of a representation $V$ of
$U_q(\widehat{sl_{n+1}})$ or $U_q(sl_n)$,
$B$ is a crystal by the obvious definitions of
$wt,\et_i,\ft_i,\varphi_i,\varepsilon_i$.
\vskip3mm

\noindent
{\bf Definition A.2.}
A morphism $\psi:B_1\ra B_2$ of two crystals is a map
$\psi:B_1\sqcup\{0\}\ra B_2\sqcup\{0\}$ satisfying
the following conditions:
\item{1.} $\psi(0)=0$.
\item{2.} $\psi$ commutes with all $\et_i$ and $\ft_i$.
\item{3.} For any $b\in B_1$ $wt\psi(b)=wt b$.
\par
\vskip3mm

\noindent
The morphism defined here is called "strict" in [K3].
A morphism is called an isomorphism if the associated map
$\psi:B_1\sqcup\{0\}\ra B_2\sqcup\{0\}$ is bijective.

For two crystals $B_1$ and $B_2$ we can define the tensor
product $B_1\ot B_2$. Here we only give the description of the action
of $\et_i$ and $\ft_i$.
$$
\eqalign{
\et_i(b_1\ot b_2)&=
\et_ib_1\ot b_2\quad \varphi_i(b_1)\geq\varepsilon_i(b_2),
\cr
&=b_1\ot \et_ib_2\quad \varphi_i(b_1)<\varepsilon_i(b_2),
\cr
\ft_i(b_1\ot b_2)&=
\ft_ib_1\ot b_2\quad \varphi_i(b_1)>\varepsilon_i(b_2),
\cr
&=b_1\ot \ft_ib_2\quad \varphi_i(b_1)\leq\varepsilon_i(b_2).
}
$$
\vskip3mm
\noindent
{\bf Remark A.3.}
\item{1.} A finite dimensional representation of $\widehat{sl_{n+1}}$
does not necessarily have a crystal base. For example,
the representation of $U_q(\widehat{sl_{3}})$ obtained from
the adjoint representation of $U_q(sl_{3})$ by the evaluation
homomorphism $U_q(\widehat{sl_{3}})\ra U_q(sl_{3})$ does not have
a crystal base. In particular, the crystal of the adjoint
representation of $U_q(sl_{3})$ cannot be extended to the normal
crystal (of [K3]) for $\widehat{sl_{3}}$.
\item{2.} Let $B_{l\La_k}$ be the crystal of the irreducible
highest weight representation of $U_q(sl_n)$ with
highest weight $l\La_k$. It is proved in [KMN2]
that $B_{l\La_k}$ can be extended to a classical crystal $B$
for $\widehat{sl_n}$ by defining a suitable action of $\et_0$ and $\ft_0$.
This crystal $B$ has the symmetry of a Dynkin diagram automorphism.
Namely, even if we permute, in a cyclic way,
the colors of the crystal graph $B$, the resulting crystal
is isomorphic to $B$.
%\bye
\vfill\eject

\magnification=\magstep1
\baselineskip 12pt

\noindent
{\fontB Reference}
\vskip3mm
\noindent
\item{[ABF]}Andrews, G., Baxter, R. and Forrester, P.,
{Eight-vertex SOS model and generalized Rogers-Ramanujan
type identities,
J. Stat. Phys.{\bf 35}, 193-266 (1984)}
\vskip2mm
\item{[B]} Butler, L.M., {Subgroup lattices and symmetric functions,
Memoir of the AMS {\bf Vol. 112} No. 539 (1994)}
\vskip2mm
\item{[D]} Dasmahapatra, S., {On the combinatorics of row and corner
transfer matrices of the $A_{n-1}^{(1)}$ restricted face models,
hep-th/9512095}
\vskip2mm
\item{[DJKMO]} Date, E, Jimbo, M., Kuniba, A., Miwa, T. and Okado, M.
{One dimensional configuration sums in vertex models
and affine Lie algebra characters, LMP {\bf 17}, 69-77 (1989)
\vskip2mm
\item{[DJO]}Date, E., Jimbo, M., Okado, M.,
{Crystal base and q vertex operators, CMP {\bf 155}, 47-69 (1993)}
\vskip2mm
\item{[JMMO]}Jimbo, M., Misra, K., Miwa, T., Okado, M.,
{Combinatorics of representation of $U_q(\hat{sl}(n))$
at $q=0$, CMP {\bf 136}, 543-566 (1991)}
\vskip2mm
\item{[K1]}Kashiwara, M.,
{On crystal bases of the q-analogue of universal
enveloping algebras, Duke Math. J. {\bf 63} No.2, 465-516 (1991)}
\vskip2mm
\item{[K2]}Kashiwara, M.,
{The crystal base and Littelmann's refined Demazure
character formula, Duke Math. J. {\bf 71} No.3, 839-858 (1993)}
\vskip2mm
\item{[K3]}Kashiwara, M.,
{Crystal bases of modified quantized
enveloping algebras, Duke Math. J. {\bf 73} No.2, 383-413 (1994)}
\vskip2mm
\item{[KMN1]}Kang, S-J., Kashiwara, M., Misra, K.,
Miwa, T., Nakashima, T., Nakayashiki, A.
{Affine crystals and vertex models,
Int. J. Mod. Phys. A {\bf 7}, Suppl. 1A, 449-484 (1992)}
\vskip2mm
\item{[KMN2]}Kang, S-J., Kashiwara, M., Misra, K.,
Miwa, T., Nakashima, T., Nakayashiki, A.
{Perfect crystals of quantum affine Lie algebras,
Duke Math. J.{\bf 68} No.3, 499-607 (1992)}
\vskip2mm
\item{[KN]} Kashiwara, M. and Nakashima, T.,{Crystal graph
for representations of the q-analogue of classical Lie algebras,
J. Alg. {\bf 165}, 295-345 (1994)}
\vskip2mm
\item{[Kir]}Kirillov, A. N.,
{Dilogarithm identities, preprint,hep-th/9408113 (1994)}
\vskip2mm
\item{[KR]} Kirillov, A. N. and Reshetikhin, N. Yu.,
{The Bethe ansatz and the combinatorics of Young tableaux,
J. Soviet Math. {\bf 41}, 925-955 (1988)}
\vskip2mm
\item{[L]} Lusztig, L.,{Green polynomials and singularities of
nilpotent classes, Adv. Math., {\bf 42}, 173-204 (1981)}
\vskip2mm
\item{[LLT]} Lascoux, A., Leclerc, B. and Thibon, J-Y.,
{Crystal graphs and $q$-analogue of weight multiplicities for the
root system $A_n$, Lett. Math. Phys. 359-374 (1995)}
\vskip2mm
\item{[LS]} Lascoux, A. and Schutzenberger, M.P.,
{Sur une conjecture de H.O. Foulkes, C.R. Acad. Sc. Paris
{\bf 288A}, 323-324 (1978)}
\vskip2mm
\item{[M]} Macdonald, I., {Symmetric functions and Hall polynomials,
second edition, Oxford Univ. Press (1995)}
\vskip2mm
\item{[NY1]} Nakayashiki, A. and Yamada, Y., {Crystallizing the spinon
basis, to appear in CMP, hep-th/9504052}
\vskip2mm
\item{[NY2]} Nakayashiki, A. and Yamada, Y., {Crystalline spinon basis
for RSOS models, to appear in IJMPA, hep-th/9505083}

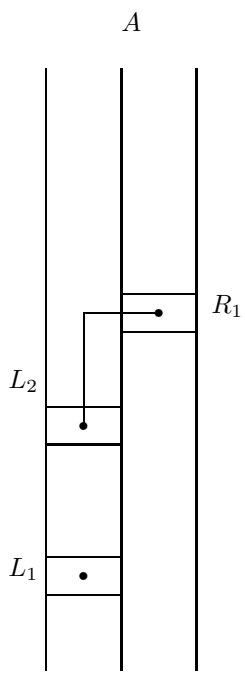
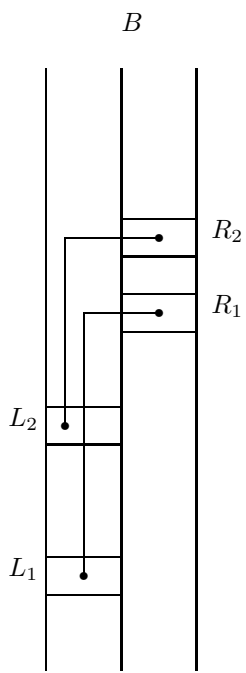
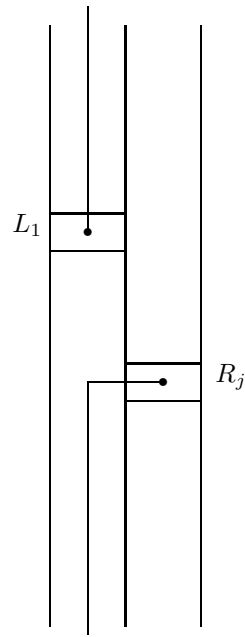
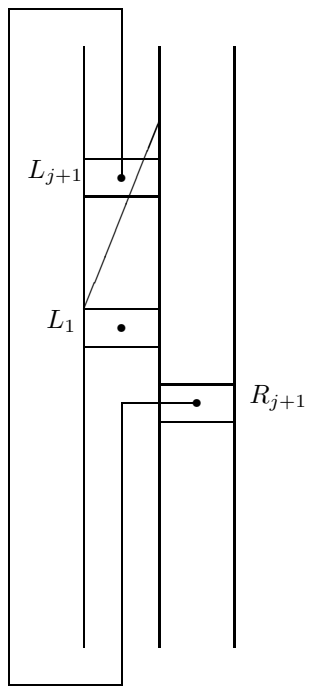
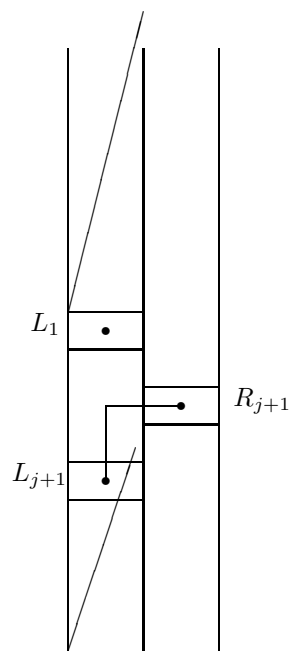
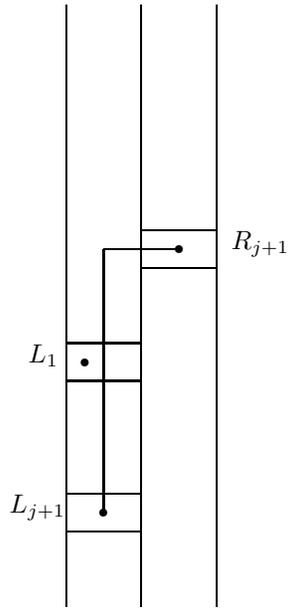
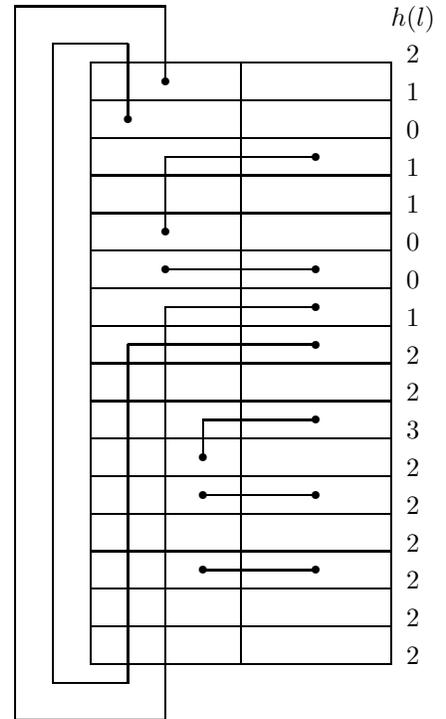
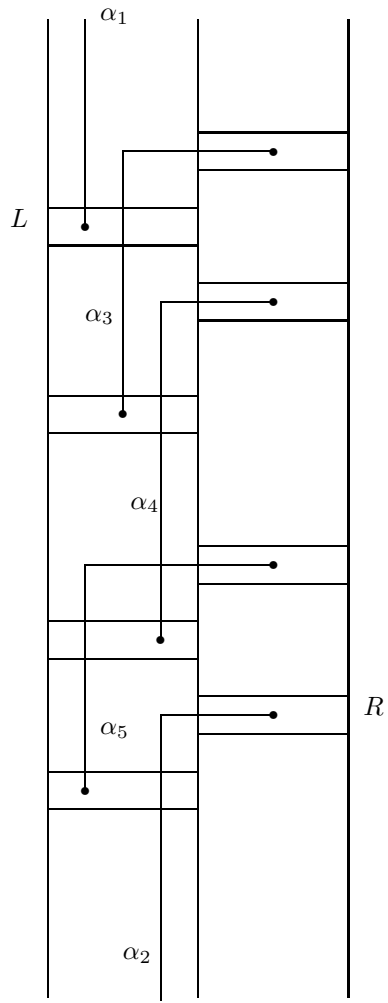
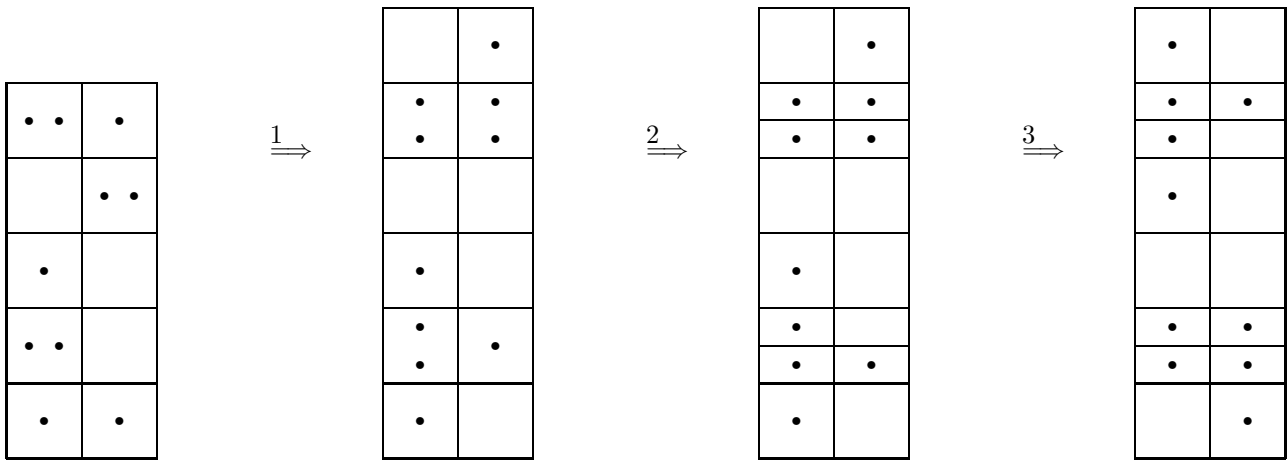
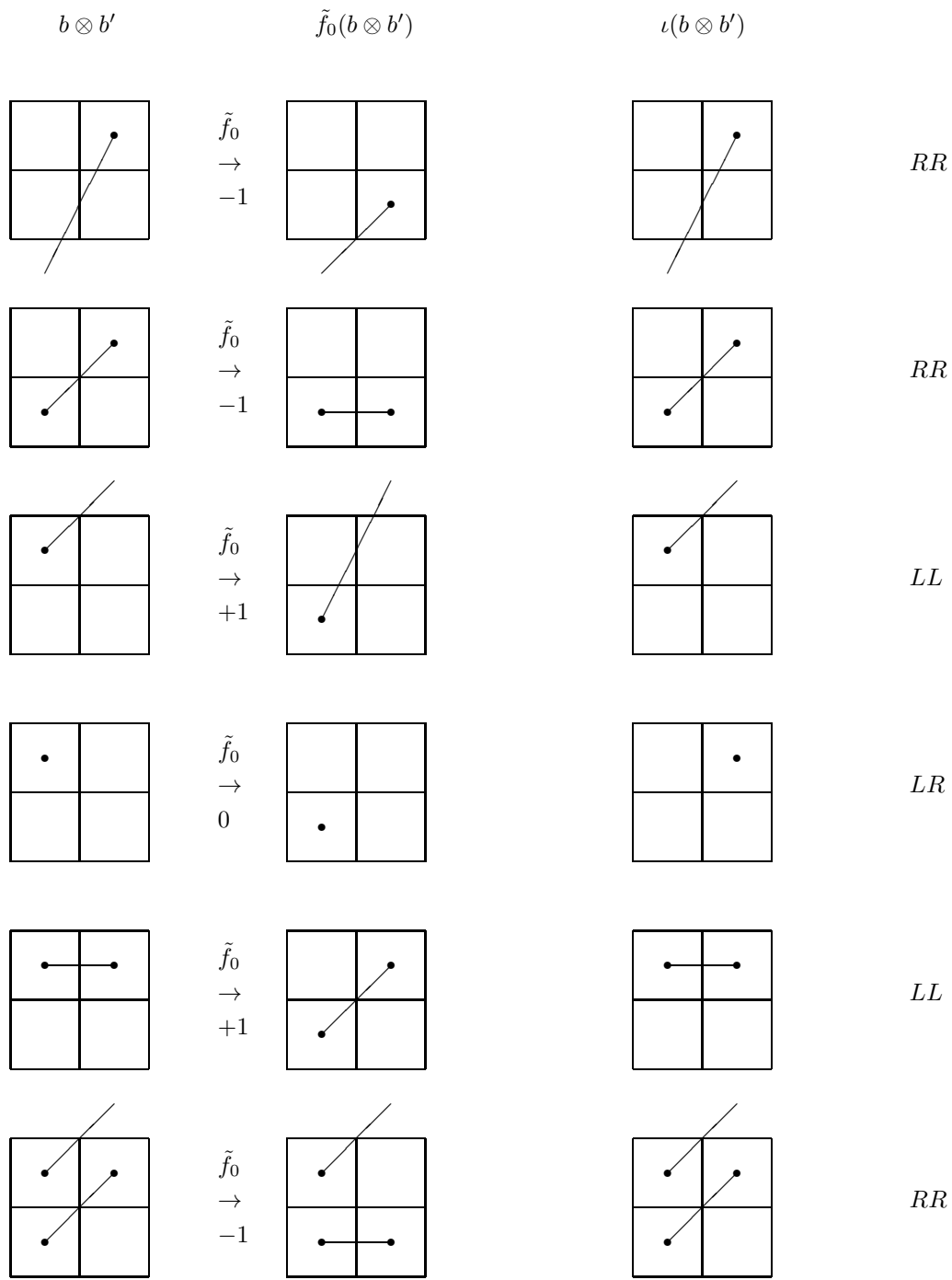
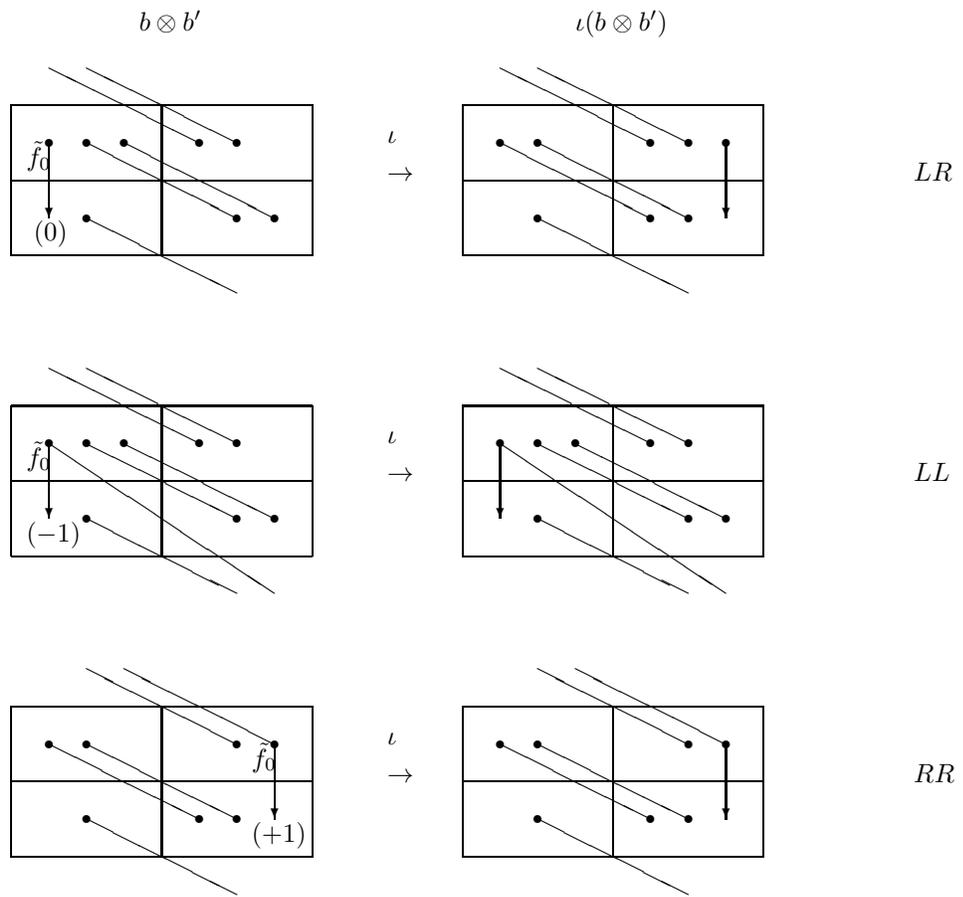
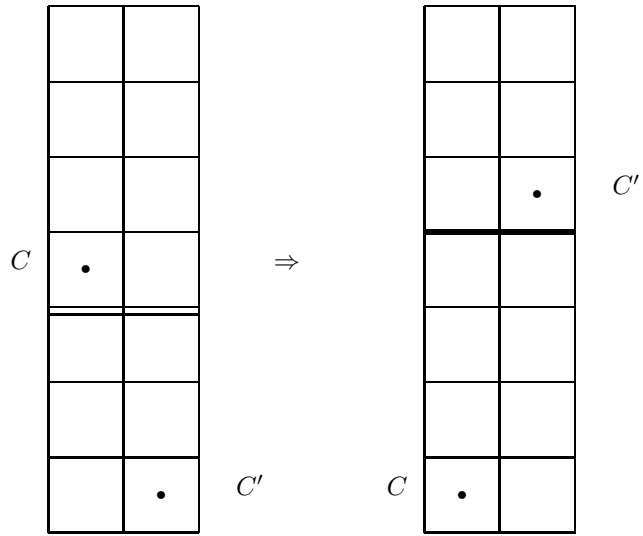
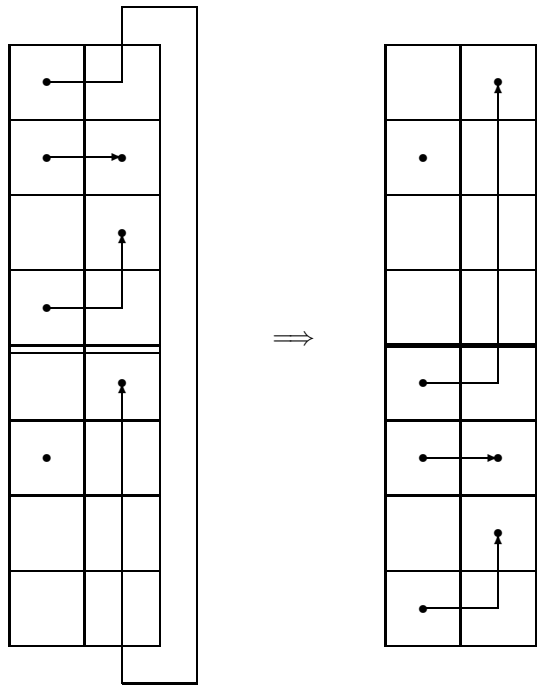
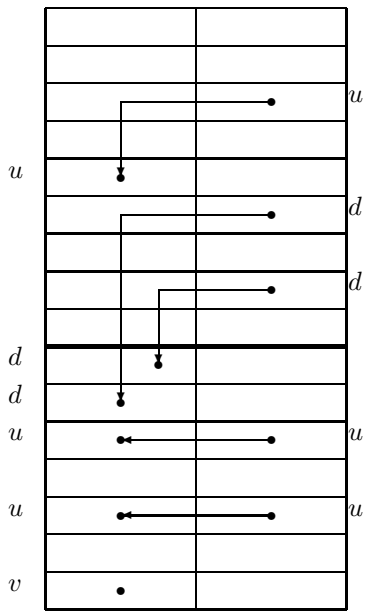
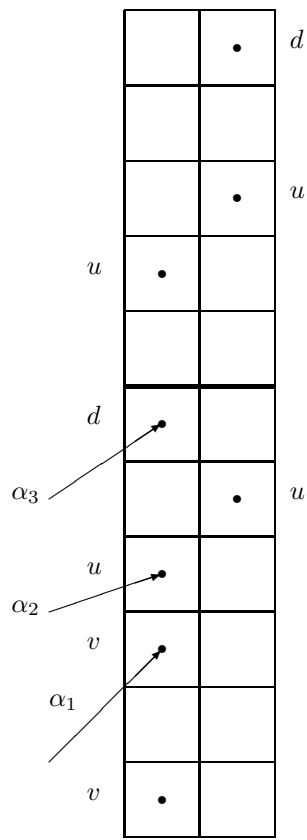
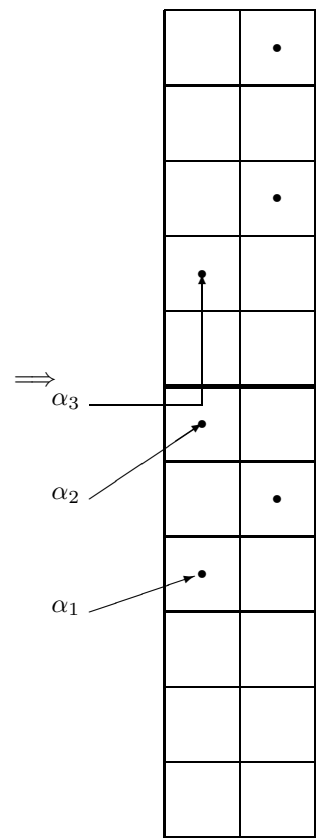
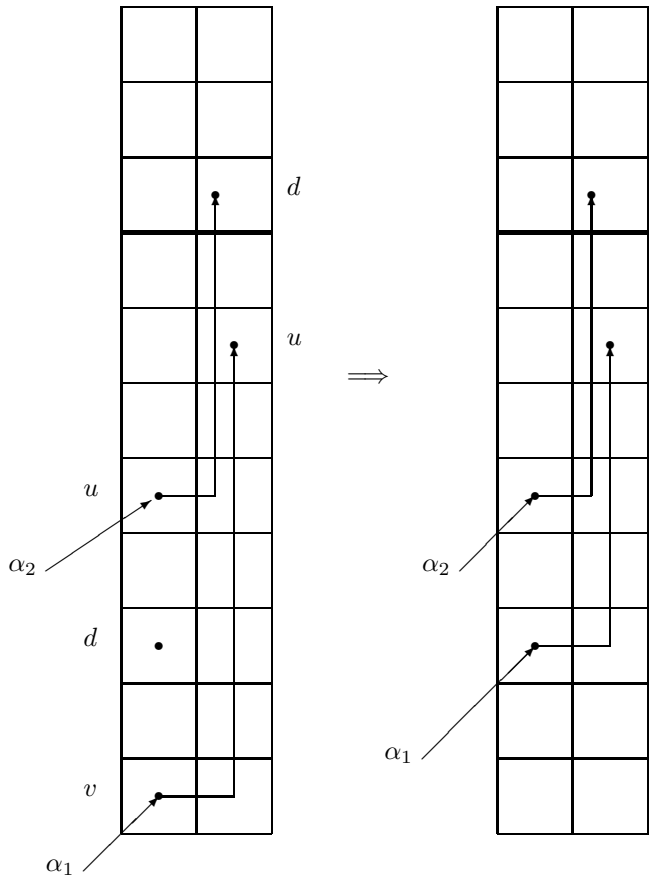
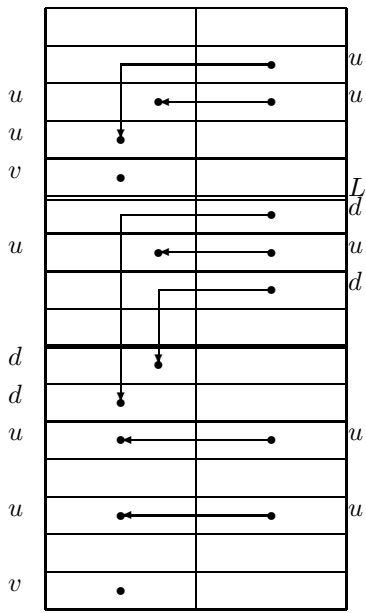
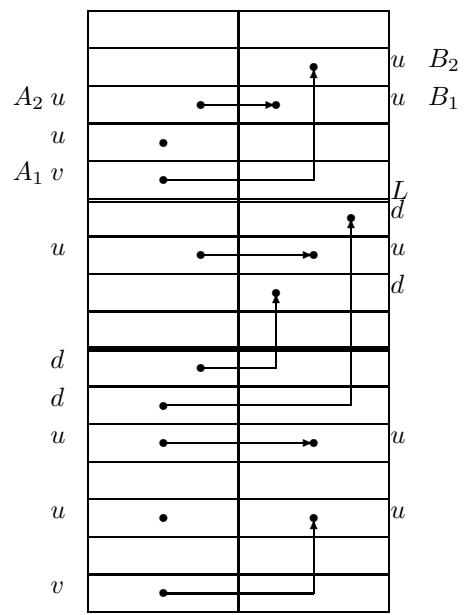
\bye